\documentclass[10.5pt,a4paper]{article}

\usepackage{graphicx,latexsym,euscript,makeidx,color,bm}
\usepackage{amsmath,amsfonts,amssymb,amsthm,thmtools,mathtools,mathrsfs,enumerate}
\usepackage[colorlinks,linkcolor=blue,citecolor=red]{hyperref}

\usepackage{geometry}
\allowdisplaybreaks[4]
  \def\cA{{\cal A}}  
    
  \def\cC{{\cal C}}  
 \def\sD{\mathscr{D}}   
\def\dbE{\mathbb{E}}    
\def\dbF{\mathbb{F}}  \def\cF{{\cal F}}  
  \def\cG{{\cal G}}  
\def\dbH{\mathbb{H}}    
\def\dbI{\mathbb{I}}

\def\dbP{\mathbb{P}}    
  \def\cQ{{\cal Q}}  
\def\dbR{\mathbb{R}}    
\def\dbS{\mathbb{S}}    
    
 \def\sU{\mathscr{U}}

             \def\hb{\hbox}
\def\ms{\medskip}              \def\ae{\hb{a.e.}}
\def\bs{\bigskip}        \def\lan{\langle}    
\def\ds{\displaystyle}   \def\ran{\rangle}    \def\tr{\hb{tr$\,$}}
         
\def\no{\noindent}          
         
\def\nn{\nonumber}         
\def\rf{\eqref}            
\def\cd{\cdot}             
\def\deq{\triangleq}     \def\({\Big (}       \def\ba{\begin{aligned}}
\def\les{\leqslant}      \def\){\Big )}       \def\ea{\end{aligned}}
\def\ges{\geqslant}      \def\[{\Big[}        \def\bel{\begin{equation}\label}
\def\ti{\tilde}          \def\]{\Big]}        \def\ee{\end{equation}}
      \def\q{\quad}        
\def\h{\widehat}         \def\qq{\qquad}      
                                              
\def\a{\alpha}  \def\G{\Gamma}      \def\Om{\Omega}  
   \def\D{\Delta}   \def\d{\delta}        
         
\def\f{\varphi}   \def\l{\lambda}        \def\e{\varepsilon}
\def\t{\tau}          

\newtheoremstyle{thry}
{}      
{}      
{\sl}   
{}      
{\bf}   
{.}     
{.5em}  
{}      
\theoremstyle{thry}

\newtheorem{theorem}{Theorem}[section]
\newtheorem{proposition}[theorem]{Proposition}

\newtheorem{lemma}[theorem]{Lemma}

\theoremstyle{definition}
\newtheorem{definition}[theorem]{Definition}

\theoremstyle{remark}
\newtheorem{remark}[theorem]{Remark}

\def\punct{}
\newtheoremstyle{dotless}{}{}{\rm}{}{\bf}{\punct}{.5em}{}
\theoremstyle{dotless}

\newtheorem{taggedassumptionx}{}
\newenvironment{taggedassumption}[1]
 {\renewcommand\thetaggedassumptionx{#1}\taggedassumptionx}
 {\endtaggedassumptionx}

\makeatletter
   
   \@addtoreset{equation}{section}
   \newcommand{\setword}[2]{%
   \phantomsection
   #1\def\@currentlabel{\unexpanded{#1}}\label{#2}%
   }
\makeatother

\begin{document}

\title{\bf  Mean-Field Stochastic Linear-Quadratic Optimal Controls:
Roles of Expectation  and Conditional Expectation Operators}

\author{Hanxiao Wang\thanks{ School of Mathematical Sciences, Shenzhen University, Shenzhen,
518060, China (Email: {\tt hxwang@szu.edu.cn}). This author is supported in part by NSFC Grants 12201424, 12426311,
Guangdong Basic and Applied Basic Research Foundation 2023A1515012104,
the Science and Technology Program of Shenzhen RCBS20231211090537064,
and Shenzhen University 2035 Program for Excellent Research Grant 2024C008.}
                           ~~~and~~
       Jiongmin Yong\thanks{Department of Mathematics, University of Central Florida,
                         Orlando, FL 32816, USA (Email: {\tt jiongmin.yong@ucf.edu}).
                          This author is supported in part by NSF Grant DMS-2305475.}}

\ms

\maketitle

\centerline{({\it Dedicated to Professor Xun Yu Zhou on the occasion of his 60th birthday})}

\bs

\no{\bf Abstract.}
This paper investigates a mean-field linear-quadratic optimal control problem where the state dynamics and cost functional incorporate both expectation and conditional expectation terms. We explicitly derive the pre-committed, na\"{\i}ve, and equilibrium solutions and establish the well-posedness of the associated Riccati equations. This reveals how the expectation and conditional expectation operators influence time-consistency.

\ms

\no{\bf Key words.}
mean-field stochastic differential equation, linear-quadratic, optimal control, time-inconsistency.

\ms

\no{\bf AMS subject classifications.}  49N10, 49N80, 60H10, 93E20.

\ms

\section{Introduction}\label{Sec:Intro}
Let $(\Om,\cF,\dbP)$ be a complete probability space on which a standard one-dimensional
Brownian motion $W(\cd)=\{W(t);0\les t<\infty\}$ is defined.
The augmented natural filtration of $W(\cd)$ is denoted by $\dbF=\{\cF_t\}_{t\ges0}$.
For any given initial pair $(t,\xi)\in\sD$ with
$$\sD=\Big\{(t,\xi)\bigm|t\in[0,T),~\xi\in L_{\cF_t}^2(\Omega;\dbR^n)\Big\},$$
consider the following controlled linear mean-field stochastic differential equation
(SDE, for short) on the finite horizon $[t,T]$:
\bel{state}\left\{\begin{aligned}
   dX(s) &=\big\{A(s)X(s)+ \bar A(s)\dbE_t[X(s)]+\ti A(s)\dbE[X(s)]+B(s)u(s) \big\}ds\\
   &\q+\big\{C(s)X(s)+ \bar C(s)\dbE_t[X(s)]+\ti C(s)\dbE[X(s)]+D(s)u(s) \big\}dW(s),\q s\in[t,T], \\
    X(t) &= \xi,
\end{aligned}\right.\ee
where $A,\bar A,\ti A,C,\bar C,\ti C:[0,T]\to\dbR^{n\times n}$, $B,D:[0,T]\to\dbR^{n\times m}$,
 called the {\it coefficients} of the {\it state equation} \rf{state}, are given deterministic functions.
The solution $X(\cd)$ of \rf{state} is called a {\it state process},
and $u(\cd)$, belonging to the space
$$ \sU[t,T] = \bigg\{\f:[t,T]\times\Om\to\dbR^{m} \bigm| \f\hb{~is $\dbF$-progressively measurable},
~\dbE\[\int^T_t|\f(s)|^2ds\]<\infty\bigg\},
$$
is called the {\it control process}.
To measure the performance of the control $u(\cd)$, we introduce the following cost functional:
\begin{align}\label{cost}
J(t,\xi;u(\cd))
&= \dbE\[\int_t^T \(\lan QX,X\ran +\lan\bar Q\dbE_t[X],\dbE_t[X]\ran
+\lan\ti Q\dbE[X],\dbE[X]\ran+\lan Ru,u\ran \) ds\nn\\
&\qq +\lan GX(T),X(T)\ran +\lan\bar G\dbE_t[X(T)],\dbE_t[X(T)]\ran
+\lan\ti G\dbE[X(T)],\dbE[X(T)]\ran  \],
\end{align}
where $G$, $\bar G$, and $\ti G$ are $n\times n$ symmetric matrices;
$Q,\bar Q,\ti Q:[0,T]\to\dbR^{n\times n}$ and $R:[0,T]\to\dbR^{m\times m}$ are deterministic symmetric matrix-valued functions.
In the Lebesgue integral on the right-hand side of \rf{cost},
we have suppressed the argument $s$, and we will do so in the sequel as long as no ambiguity arises.
With the state equation \rf{state} and   cost functional \rf{cost},
the problem can be stated as follows:

\ms

{\bf Problem (MF-SLQ).} For any given initial pair $(t,\xi)\in\sD$,
find a control $\bar u(\cd)\in \sU[t,T]$ such that
\bel{inf J}J(t,\xi;\bar u(\cd))=\inf_{u(\cd)\in\sU[t,T]} J(t,\xi;u(\cd)).\ee

Any $\bar u(\cd)\in\sU[t,T]$ satisfying \rf{inf J} is called an {\it optimal control}
of Problem (MF-SLQ) for the initial pair $(t,\xi)$;
the corresponding state process $\bar X(\cd)\equiv X(\cd\,;t,\xi,\bar u(\cd))$ is called an
{\it optimal state process}; and
$(\bar X(\cd),\bar u(\cd))$ is called an  {\it optimal pair}.
Compared with classical stochastic LQ optimal control problems,
the main feature of our Problem (MF-SLQ) is that the state dynamics and cost functional incorporate both expectation and conditional expectation terms of the state process. Such a type of problems is usually referred to as a {\it mean-field stochastic linear-quadratic (SLQ) optimal control problem.}

The mean-field SLQ optimal control problem was initially studied by Yong \cite{Yong2013},
in which the state dynamics and cost functional depend on the expectation of the state and control processes.
Such a type of control problems is mainly motivated by the famous {\it mean-field game} or {\it mean-field control} theory,
which was independently introduced by Huang--Malham\'{e}--Caines \cite{Huang2006} and Lasry--Lions \cite{Lasry2007}
in the early 2000s. The mean-field term, that is the expectation term,
comes from the limit of the average state in large population models.
We refer the reader to \cite{Pham2016,Sun2017,Basei2019,Sun2021,Huang2024,Moon-Basar2024,Li-Li-Xu2025}
for the newest developments of mean-field SLQ optimal control problems along this line.

Another motivation for studying mean-field SLQ control problems is
the so-called {\it dynamic mean-variance model} (see Basak--Chabakauri \cite{Basak-Chabakauri2010}, for example),
in which the cost functional depends on the conditional expectation of the terminal state in a nonlinear way.
In this type of problems, the optimal control is not expected to be time-consistent.
In other words, an optimal control selected at a given initial pair $(t,\xi)\in\sD$  will not stay optimal thereafter.
Thus, instead of looking for optimal controls, one usually hopes to find an equilibrium strategy.

The earliest mathematical consideration of time-inconsistent problems
was given by Strotz \cite{Strotz1955}, followed by Pollak \cite{Pollak1968}, and the recent works of   Ekeland--Lazrak \cite{Ekeland2010}, Yong \cite{Yong2012,Yong2014}, Bj\"{o}rk--Khapko--Murgoci \cite{Bjork-Khapko-Murgoci2017},
 He--Jiang \cite{He-Jiang2021},
Hern\'{a}ndez--Possama\"{i} \cite{Hernandez-2020}, and Wang--Yong--Zhou \cite{Wang2024}, to mention a few.
 For mean-field SLQ control problems, the related study from a time-inconsistent viewpoint can be found in Basak--Chabakauri \cite{Basak-Chabakauri2010}, Hu--Jin--Zhou \cite{Hu-Jin-Zhou2012,Hu-Jin-Zhou2017},
 Bj\"ork--Murgoci--Zhou \cite{Bjork-Murgoci-Zhou2014}, Yong \cite{Yong2017}, Wang \cite{Wang2020},
   L\"{u}--Ma--Wang \cite{Lv2025}, etc.

As noted above, there are two distinct formulations of mean-field SLQ optimal control problems in the literature:
One uses expectation terms \cite{Yong2013, Pham2016,Sun2017,Basei2019,Sun2021,Huang2024,Moon-Basar2024,Li-Li-Xu2025},
while the other employs conditional expectation terms \cite{Basak-Chabakauri2010,Hu-Jin-Zhou2012,Bjork-Murgoci-Zhou2014,Hu-Jin-Zhou2017,Yong2017,Wang2020,Lv2025},
and these have long been studied independently.
In this paper, we incorporate  the expectation and conditional expectation terms together into one controlled system, and hope to provide a relatively  uniform approach to these two types of problems.
The {\it pre-committed}, {\it na\"{\i}ve}, and {\it equilibrium} solutions, initially introduced
by Strotz \cite{Strotz1955}, will be all derived explicitly,
in terms of the solutions to the associated Riccati equations. Consequently,
 the roles of  expectation and conditional expectation operators can be clarified.

\ms

The main results of this paper can be  summarized as follows.

\begin{enumerate}[(i)]
\item The optimal control $u^*(\cd)$ of Problem (MF-SLQ) at $(t,\xi)\in\sD$ can be explicitly given by the following form:
\bel{u*}
u^*(\cd)=\Psi^*(\cd)X^*(\cd)+\bar\Psi^*(\cd)\dbE_t[X^*(\cd)]+\ti \Psi^*(\cd)\dbE[X^*(\cd)],\ee
where $(\Psi^*(\cd),\bar\Psi^*(\cd),\ti \Psi^*(\cd))$ are determined by the unique solution of the associated three Riccati equations \rf{RE-pre-1}--\rf{RE-pre-3},
and $X^*(\cd)$ is the corresponding optimal state process.

We can see that the {\it optimal feedback strategy} $(\Psi^*(\cd),\bar\Psi^*(\cd),\ti \Psi^*(\cd))$ is time-consistent,
because it is independent of the initial time $t$ (which can be regarded as the current time), while the {\it optimal feedback operator} $\mathbf{\Psi_t^*}[\,\cd\,]$,
defined by
\begin{align}
[\mathbf{\Psi_t^*}X](s)&:=[\Psi^*(s)+\bar\Psi^*(s)\dbE_t+\ti \Psi^*(s)\dbE][X(s)]\nn\\
&=\Psi^*(s)X(s)+\bar\Psi^*(s)\dbE_t[X(s)]+\ti \Psi^*(s)\dbE[X(s)],\q s\in[t,T],\nn\\
&\qq\qq\forall X(\cd)\in L_\dbF^2(\Omega;C([t,T];\dbR^n)),\label{pre-opti-operator}
\end{align}
turns out to be time-inconsistent, mainly due to the fact that it depends on $\dbE_t[\,\cd\,]$. Thus, Problem (MF-SLQ) is time-inconsistent.
Consequently, the explicit representation \rf{u*} only guarantees $u^*(\cd)$ to be optimal at the current moment $t$.
Because of this, we call $u^*(\cd)$ given by \rf{u*} a {\it pre-committed control} over $[t,T]$,
and the corresponding state process $X^*(\cd)$ a {\it pre-committed state} over $[t,T]$.

\item
Since the optimal feedback operator $\mathbf{\Psi_t^*}[\,\cd\,]$ depends on $t$ through the conditional expectation operator $\dbE_t[\,\cd\,]$,
Problem (MF-SLQ) is time-inconsistent.
Motivated by \cite{Strotz1955,Chen-Zhou2024}, a {\it na\"{\i}ve feedback operator} is introduced.
It is constructed as the limit of a sequence of optimal feedback operators over every small time interval.
Precisely, the na\"{\i}ve feedback operator $\mathbf{\h\Psi}[\,\cd\,]$ is given by
\begin{align}
[\mathbf{\h\Psi}X](s)&=[\h\Psi_1(s)+\h\Psi_2(s)\dbE][X(s)]\nn\\
&=[\Psi^*(s)+\bar\Psi^*(s)]X(s)+\ti \Psi^*(s)\dbE[X(s)],\q s\in[0,T],\nn\\
&\qq\qq\forall X(\cd)\in L_\dbF^2(\Omega;C([0,T];\dbR^n)).\label{naive-operator}
\end{align}
By comparing \rf{pre-opti-operator} and \rf{naive-operator}, we  see that
$\mathbf{\h\Psi}[\,\cd\,]$ can be deduced from the pre-committed solution $\mathbf{\Psi_0^*}[\,\cd\,]$ over $[0,T]$
 by replacing $\dbE_0[\,\cd\,]$ with $\dbE_s[\,\cd\,]$.

\item The {\it equilibrium feedback operator} of Problem (MF-SLQ) is constructed by the multi-person
differential game method introduced by Yong \cite{Yong2017}, which admits the following form:
\begin{align}
[\mathbf{\Psi^\dag}X](s)&:=[\Psi_1^\dag(s)+ \Psi_2^\dag(s)\dbE][X(s)]\nn\\
&=[\Psi^\dag(s)+\bar\Psi^\dag(s)]X(s)+\ti \Psi^\dag(s)\dbE[X(s)],\q s\in[0,T],\nn\\
&\qq\qq\forall X(\cd)\in L_\dbF^2(\Omega;C([0,T];\dbR^n)),\label{equilibrium-operator}
\end{align}
with  $(\Psi^\dag(\cd),\bar\Psi^\dag(\cd),\ti \Psi^\dag(\cd))$ being given by
the unique solution of three equilibrium Riccati equations \rf{ER-1}--\rf{ER-3}.
Comparing with Yong \cite{Yong2017}, we see that \rf{equilibrium-operator} additionally depends on
the expectation of the state process, and some unreasonable assumptions
(i.e., (H4) in \cite{Yong2017}) are removed by
fixing a crucial gap.
\end{enumerate}

From the above, we see that only the conditional expectation operator $\dbE_t[\,\cd\,]$ affects the time-consistency of Problem (MF-SLQ).
In the pre-committed solution $\mathbf{\Psi_0^*}[\,\cd\,]$ over $[0,T]$, the feedback strategy $\bar\Psi^*(\cd)$ acts on $\dbE[X(\cd)]$. This indicates that the conditional expectation operator functions similarly to the standard expectation operator within this solution.
Conversely, in both the na\"{\i}ve and equilibrium solutions, the respective feedback strategies $\bar\Psi^*(\cd)$ and $\bar\Psi^\dag(\cd)$ are applied directly to the state process $X(\cd)$. Here, the conditional expectation term functions similarly to the state term itself in conventional problems.
However, $\bar\Psi^*(\cd)$ and $\bar\Psi^\dag(\cd)$ are constructed by different Riccati equations.

The rest of this paper is organized as follows.
In Section \ref{sec:main-results}, we state the main results of our paper, with some explanations.
In Section \ref{sec:nai}, we give the convergence analysis of the  na\"{\i}ve solution.
The equilibrium solution is derived in Section \ref{sec:equilib},
and the verification theorem of equilibrium solutions is presented in Section \ref{sec:verifi}.
Some lengthy proofs are given in Appendix.

\section{The Main Results}\label{sec:main-results}

Let $\dbR^{n\times m}$ be the Euclidean space consisting of $n\times m$ real matrices,
endowed with the Frobenius inner product $\lan M,N\ran\deq\tr[M^\top N]$,
where $M^\top$ and $\tr(M)$ stand for the transpose  and  the trace of  $M$, respectively.
Let $\dbS^n$  be the subspace of $\dbR^{n\times n}$ consisting of symmetric matrices
and $\dbS^n_+$ the subset of $\dbS^n$ consisting of positive semi-definite matrices.
For any Euclidean space $\dbH$ (which could be $\dbR^n$, $\dbR^{n\times m}$, $\dbS^n$, etc.),
we introduce the following spaces:
\begin{align*}
C([0,T];\dbH):
   &\hb{~~the space of $\dbH$-valued, continuous functions on $[0,T]$}; \\
L^\infty([0,T];\dbH):
   &\hb{~~the space of $\dbH$-valued,  essentially bounded functions  on $[0,T]$};\\
L^2_{\mathcal{F}_t}(\Om;\dbH):
   &\hb{~~the space of $\mathcal{F}_t$-measurable, $\dbH$-valued, square-integrable random variables}; \\
L_\dbF^2([0,T];\dbH):
   &\hb{~~the space of $\dbF$-progressively measurable, $\dbH$-valued processes} \\
   &\hb{~~$\f:[0,T]\times\Om\to\dbH$ with $\ds\dbE\[\int_0^T|\f(s)|^2ds\]<\infty$};\\
L_\dbF^2(\Om;C([0,T];\dbH)):
   &\hb{~~the space of $\dbF$-adapted, continuous, $\dbH$-valued processes} \\
   &\hb{~~$\f:[0,T]\times\Om\to\dbH$ with $\dbE\[\ds\sup_{s\in[0,T]}|\f(s)|^2\]<\infty$}.
\end{align*}
For $M,N\in\dbS^n$, we use the notation $M\ges N$ (respectively, $M>N$)
to indicate that $M-N$ is positive semi-definite (respectively, positive definite).
For any $\dbS^n$-valued measurable function $F(\cd)$ on $[0,T]$, we denote
$$\left\{\begin{aligned}
  F(\cd)\ges 0\q &\Longleftrightarrow\q F(s)\ges 0,\q \ae~s\in[0,T];\\
 F(\cd)> 0\q &\Longleftrightarrow\q F(s)> 0,\q \ae~s\in[0,T];\\
 F(\cd)\gg 0\q &\Longleftrightarrow\q F(s)\ges \d I_n,\q \ae~s\in[0,T],~\hbox{for some } \d>0.
\end{aligned}\right. $$

\ms
To guarantee that Problem (MF-SLQ) is well-posed,
we impose the following assumptions for the state equation \rf{state}
and cost functional \rf{cost}.

\begin{taggedassumption}{(H1)}\label{ass:A1}
The coefficients of state equation \rf{state} satisfy
$$
A(\cd),\bar A(\cd),\ti A(\cd), C(\cd),\bar C(\cd),\ti C(\cd) \in L^\infty([0,T];\dbR^{n\times n}),
\q  B(\cd), D(\cd)\in L^\infty([0,T];\dbR^{n\times m}).
$$
\end{taggedassumption}

\begin{taggedassumption}{(H2)}\label{ass:A2}
The weighting coefficients in the quadratic functional \rf{cost} satisfy:
$G,\bar G,\ti G\in\dbS_+^n$,
$Q(\cd),\bar Q(\cd),\ti Q(\cd)\in L^\infty([0,T];\dbS^n)$, and $R(\cd)\in L^\infty([0,T];\dbS^m)$ with
 $$
 Q(\cd),\bar Q(\cd),\ti Q(\cd)\ges 0,\qq
 R(\cd)\gg0.
 $$
\end{taggedassumption}

By \cite[Lemma 2.1]{Yong2013}, we have the following result.

\begin{lemma}\label{lem:state-well}
Let {\rm\ref{ass:A1}} hold. Then for any $(t,\xi)\in\sD$ and $u(\cd)\in\sU[t,T]$,
the state equation \rf{state} admits a unique solution $X(\cd)\in L_\dbF^2(\Om;C([0,T];\dbR^n))$.
Moreover, there exists a constant $K>0$, which is independent of  $(t,\xi)$ and $u(\cd)\in\sU[t,T]$, such that
$$
\dbE\[\sup_{s\in[t,T]}|X(s)|^2\]\les K\dbE\[|\xi|^2+\int_t^T|u(s)|^2ds\].
$$
\end{lemma}

Note that under \ref{ass:A1}, the unique solution $X(\cd)$ of \rf{state} belongs to $L_\dbF^2(\Om;C([0,T];\dbR^n))$.
In addition, if the assumption \ref{ass:A2} holds,
then the random variables on the right-hand side of  \rf{cost} are integrable and Problem (MF-SLQ) is well-posed.

The following lemma will play an important role in the construction of equilibrium solutions.
Consider the equation
\bel{state-cX}\left\{\begin{aligned}
   dX(s) &=\big\{\cA X+ \bar \cA\dbE_t[X]+\ti\cA\dbE[X] \big\}ds\\
   &\qq+\big\{\cC X+ \bar \cC\dbE_t[X]+\ti\cC\dbE[X]\big\}dW(s),\q s\in[t,T], \\
    X(t) &= \xi,
\end{aligned}\right.\ee
and the functional
\begin{align}
I(X(\cd))&=\dbE\[\int_t^T\(\lan \cQ_1 X,X\ran+\lan \cQ_2\dbE_t[X],\dbE_t[X]\ran+\lan \cQ_3\dbE_\t[X],\dbE_\t[X]\ran+\lan \cQ_4\dbE[X],\dbE [X]\ran \)ds\nn\\
&\qq +\lan \cG_1 X(T),X(T)\ran+\lan \cG_3\dbE_\t[X(T)],\dbE_\t[X(T)]\ran+\lan \cG_4\dbE[X(T)],\dbE[X(T)]\ran\],\label{I(X)}
\end{align}
where $\t\in[0,t]$ is a given time.

\begin{lemma}\label{lem:rep}
Suppose that
$\cA(\cd),\bar \cA(\cd),\ti \cA(\cd),\cC(\cd),\bar \cC(\cd),\ti \cC(\cd)\in L^\infty([0,T];\dbR^{n\times n})$,
$\cG_1,\cG_3,\cG_4\in\dbS_+^n$, and
$\cQ_1(\cd),\cQ_2(\cd),\cQ_3(\cd),\cQ_4(\cd)\in L^\infty([0,T];\dbS_+^n)$.
Then the following system of Lyapunov equations admits a unique solution
$(\check\G(\cd),\G(\cd),\bar\G(\cd), \ti\G(\cd))\in [C([0,T];\dbS_+^n)]^4$:
\bel{lem:rep1}\left\{\begin{aligned}
   &\dot{\check\G}+\check\G\cA +\cA^\top \check\G+ \cC^\top \check\G\cC+\cQ_1=0, \\
   &\check\G(T)=\cG_1,
\end{aligned}\right.\ee
\bel{lem:rep2}\left\{\begin{aligned}
   &\dot{\G}+\G(\cA+\bar\cA) +(\cA+\bar\cA)^\top \G+ (\cC+\bar\cC)^\top \check\G(\cC+\bar\cC)+\cQ_1+\cQ_2=0, \\
   &\G(T)=\cG_1,
\end{aligned}\right.\ee
\bel{lem:rep3}\left\{\begin{aligned}
   &\dot{\bar\G}+\bar\G(\cA+\bar\cA) +(\cA+\bar\cA)^\top \bar\G+ (\cC+\bar\cC)^\top \check\G(\cC+\bar\cC)+\cQ_1+\cQ_2+\cQ_3=0, \\
   &\bar\G(T)=\cG_1+\cG_3,
\end{aligned}\right.\ee
\bel{lem:rep4}\left\{\begin{aligned}
   &\dot{\ti\G}+\ti\G(\cA+\bar\cA+\ti\cA) +(\cA+\bar\cA+\ti\cA)^\top \ti\G+ (\cC+\bar\cC+\ti\cC)^\top \check\G(\cC+\bar\cC+\ti\cC)\\
   &\q+\cQ_1+\cQ_2+\cQ_3+\cQ_4=0, \\
   &\ti\G(T)=\cG_1+\cG_3+\cG_4.
\end{aligned}\right.\ee
Moreover, we have the following representation of the functional \rf{I(X)}:
\begin{align}
I(X(\cd))&=\dbE\[\big\lan \G(t)(X(t)-\dbE_\t[X(t)]),\,(X(t)-\dbE_\t[X(t)])\big\ran+\big\lan \bar\G(t)(\dbE_\t[X(t)]-\dbE[X(t)]),\nn\\
&\qq\q(\dbE_\t[X(t)]-\dbE[X(t)])\big\ran+\big\lan \ti\G(t)\dbE[X(t)],\,\dbE[X(t)]\big\ran\].\label{lem:rep5}
\end{align}
\end{lemma}

\begin{remark}
Lemma \ref{lem:rep} is a modification of Yong \cite[Lemma 2.4]{Yong2017}.
To address the expectation term $\dbE[X(\cd)]$, we additionally introduce the fourth Lyapunov equation \rf{lem:rep4}.
The proof of Lemma \ref{lem:rep} can be found in Appendix (A1).
\end{remark}

\subsection{Pre-committed Solutions}
In this subsection, we shall consider the pre-committed solution.

\begin{proposition}\label{Prop:pre-solution1}
Let {\rm\ref{ass:A1}--\ref{ass:A2}} hold. Then $u^*(\cd)\in\sU[t,T]$ is an optimal control for the initial pair $(t,\xi)\in\sD$ if and only if it satisfies
\bel{Prop:pre-solution1-main1}
R(s)u^*(s)+B(s)^\top Y(s)+D(s)^\top Z(s)=0,\q s\in[t,T],
\ee
where $(Y(\cd),Z(\cd))$ is the unique solution to the mean-field BSDE:
\bel{Prop:pre-solution1-main2}\left\{\begin{aligned}
   dY(s) &=-\big\{A^\top Y+ \bar A^\top \dbE_t[Y]+\ti A^\top\dbE[Y]+C^\top Z+ \bar C^\top \dbE_t[Z]+\ti C^\top\dbE[Z]\\
   &\qq+QX^*+ \bar Q \dbE_t[X^*]+\ti Q\dbE[X^*]\big\}ds+ZdW(s),\q s\in[t,T], \\
    Y(T) &= GX^*(T)+ \bar G \dbE_t[X^*(T)]+\ti G\dbE[X^*(T)],
\end{aligned}\right.\ee
with $X^*(\cd)$ being the state process of \rf{state} corresponding to $u^*(\cd)$, that is
\bel{Prop:pre-solution1-main3}\left\{\begin{aligned}
   dX^*(s) &=\big\{AX^*+ \bar A\dbE_t[X^*]+\ti A\dbE[X^*]+Bu^* \big\}ds\\
   &\qq+\big\{CX^*+ \bar C\dbE_t[X^*]+\ti C\dbE[X^*]+Du^* \big\}dW(s), \q s\in[t,T],\\
    X^*(t) &= \xi.
\end{aligned}\right.\ee
\end{proposition}

The proof of Proposition \ref{Prop:pre-solution1} is standard.
For completeness, we sketch it in Appendix (A2).
Note that under \ref{ass:A2}, the functional $u(\cd)\mapsto J(t,\xi;u(\cd))$ is uniformly convex.
Then by following Sun--Yong \cite{Sun-Yong2020}, it is easy to show that
for any $(t,\xi)\in\sD$, Problem (MF-SLQ) admits a unique optimal control.

To decouple the mean-field forward-backward system \rf{Prop:pre-solution1-main1}--\rf{Prop:pre-solution1-main3},
we introduce the following Riccati equations:
\bel{RE-pre-1}\left\{\begin{aligned}
   &\dot{P}+PA+A^\top P+C^\top PC+Q\\
   &\q-(PB+C^\top PD)(R+D^\top PD)^{-1}(B^\top P+D^\top PC)=0, \\
   &P(T)=G,
\end{aligned}\right.\ee
\bel{RE-pre-2}\left\{\begin{aligned}
   &\dot{\Pi}+\Pi(A+\bar A)+(A+\bar A)^\top \Pi+(C+\bar C)^\top P(C+\bar C)+Q+\bar Q\\
   &\q-[\Pi B+(C+\bar C)^\top PD](R+D^\top PD)^{-1}[B^\top \Pi+D^\top P(C+\bar C)]=0,\\
   &\Pi(T)=G+\bar G,
\end{aligned}\right.\ee
\bel{RE-pre-3}\left\{\begin{aligned}
  &\dot{\Phi}+\Phi(A+\bar A+\ti A)+(A+\bar A+\ti A)^\top \Phi+(C+\bar C+\ti C)^\top P(C+\bar C+\ti C)+Q+\bar Q+\ti Q\\
   &\q-[\Phi B+(C+\bar C+\ti C)^\top PD](R+D^\top PD)^{-1}[B^\top \Phi+D^\top P(C+\bar C+\ti C)]=0,\\
   &\Phi(T)=G+\bar G+\ti G.
\end{aligned}\right.\ee

\begin{theorem}\label{Prop:pre-solution-clsoed}
Let {\rm\ref{ass:A1}--\ref{ass:A2}}  hold. Then the system of Riccati equations \rf{RE-pre-1}--\rf{RE-pre-3}
admits a unique solution $(P(\cd),\Pi(\cd),\Phi(\cd))\in [C([0,T];\dbS_+^n)]^3$.
Moreover, the unique optimal control $u^*(\cd)$ for the initial pair $(t,\xi)\in\sD$ admits the following closed-loop representation:
\begin{align}
u^*(s)&=-\big[R(s)+D(s)^\top P(s)D(s)\big]^{-1}\Big\{\big[B(s)^\top P(s)+D(s)^\top P(s)C(s)\big]X^*(s)\nn\\
&\qq +\big[B(s)^\top \big(\Pi(s)-P(s)\big)+D(s)^\top P(s)\bar C(s)\big]\dbE_t[X^*(s)]\nn\\
&\qq +\big[B(s)^\top \big(\Phi(s)-\Pi(s)\big)+D(s)^\top P(s)\ti C(s)\big]\dbE[X^*(s)]\Big\}\nn\\
&=: \Psi^*(s)X^*(s)+\bar\Psi^*(s)\dbE_t[X^*(s)]+\ti \Psi^*(s)\dbE[X^*(s)],  \q s\in[t,T],\label{Prop:pre-solution-clsoed-main1}
\end{align}
with $X^*(\cd)$ being the unique solution of the closed-loop system:
\bel{state*}\left\{\begin{aligned}
   dX^*(s) &=\big\{(A+B\Psi^*)X^*+ (\bar A+B\bar\Psi^*)\dbE_t[X^*]+ (\ti A+B\ti\Psi^*)\dbE[X^*]\big\}ds\\
   &\q+\big\{(C+D\Psi^*)X^*+ (\bar C+D\bar \Psi^*)\dbE_t[X^*]\\
   &\q+ (\ti C+D\ti\Psi^*)\dbE[X^*] \big\}dW(s),\q s\in[t,T], \\
    X^*(t) &= \xi.
\end{aligned}\right.\ee
\end{theorem}

The proof of Theorem \ref{Prop:pre-solution-clsoed} will be given in Appendix (A3).

\begin{remark}
Note that the system of Riccati equations \rf{RE-pre-1}--\rf{RE-pre-3} is independent of the current time $t$.
Thus, the optimal feedback strategy  $(\Psi^*(\cd),\bar\Psi^*(\cd),\ti \Psi^*(\cd))$ is time-consistent.
However, the optimal control is not a direct outcome of the  optimal feedback strategy,
but an outcome of the optimal feedback operator $\mathbf{\Psi_t^*}$,
which is defined by
\begin{align}
[\mathbf{\Psi_t^*}X](s)&:=[\Psi^*(s)+\bar\Psi^*(s)\dbE_t+\ti \Psi^*(s)\dbE][X(s)]\nn\\
&=\Psi^*(s)X(s)+\bar\Psi^*(s)\dbE_t[X(s)]+\ti \Psi^*(s)\dbE[X(s)],\q s\in[t,T],\nn\\
&\qq\qq\forall X(\cd)\in L_\dbF^2(\Omega;C([t,T];\dbR^n)).\label{def:pre-oper}
\end{align}
Since the conditional expectation operator $\dbE_t[\,\cd\,]$ depends on $t$,
Problem (MF-SLQ) is time-inconsistent.
\end{remark}

\begin{remark}
The pre-committed solution over $[t,T]$ optimizes the cost functional solely at time $t$. In this framework, the controller adheres to the resulting policy thereafter---despite recognizing it may cease to be optimal at future times via a ``commitment device" if accessible and necessary. Thus, the unique {\it pre-committed solution} of Problem (MF-SLQ) over $[t,T]$ is given by \rf{Prop:pre-solution-clsoed-main1}.
The {\it pre-committed feedback operator} over $[t,T]$ is given by $\mathbf{\Psi^*_t}$, which is defined by \rf{def:pre-oper}.
\end{remark}

\subsection{Na\"{\i}ve Solutions}

The na\"{\i}ve controller is unaware of time inconsistency. At any given initial pair $(t,\xi)\in\sD$, he seeks an ``optimal" strategy solely for that instant, oblivious to its inevitable future abandonment. Consequently, his strategies perpetually change, and the actually implemented strategy emerges as the limiting outcome of a sequence of momentary ``optimal strategies".

\begin{definition}\label{naive-stra}
Let
$\D_N=\{0=t_0<t_1<...<t_{N-1}<t_N=T\big\}$ be a partition of $[0,T]$ with
$\|\D_N\|=\sup_{1\les i\les N} |t_i-t_{i-1}|$.
We call $\h u_{\D_N}(\cd)$ a {\it na\"{\i}ve control} associated with $\D_N$,
if $\h u_{\D_N}(\cd)|_{[t_i,t_{i+1})}$ coincides with the optimal control
associated with the cost functional $J(t_i,\h X_{\D_N}(t_{i});u(\cd))$ over every $[t_i,t_{i+1})$, where
$$
\left\{\begin{aligned}
   d\h X_{\D_N}(s) &=\big\{A\h X_{\D_N}+ \bar A\dbE_{t_i}[\h X_{\D_N}]+\ti A\dbE[\h X_{\D_N}]+B\h u_{\D_N} \big\}ds\\
   &\q+\big\{C\h X_{\D_N}+ \bar C\dbE_{t_i}[\h X_{\D_N}]+\ti C\dbE[\h X_{\D_N}]+D\h u_{\D_N} \big\}dW(s), \\
   &\qq\qq \qq s\in[t_i,t_{i+1}),\q i=0,1,...,N-2;\\
    d\h X_{\D_N}(s) &=\big\{A\h X_{\D_N}+ \bar A\dbE_{t_{N-1}}[\h X_{\D_N}]+\ti A\dbE[\h X_{\D_N}]+B\h u_{\D_N} \big\}ds\\
   &\q+\big\{C\h X_{\D_N}+ \bar C\dbE_{t_{N-1}}[\h X_{\D_N}]+\ti C\dbE[\h X_{\D_N}]+D\h u_{\D_N} \big\}dW(s), \\
   &\qq\qq \qq s\in[t_{N-1},t_N],\\
   \h X_{\D_N}(0) &= \xi.
\end{aligned}\right.
$$
Let $(\h\Psi_1(\cd),\h\Psi_2(\cd))\in L^\infty([0,T];\dbR^{m\times n})\times L^\infty([0,T];\dbR^{m\times n})$, $\h X(\cd)$ be the state process of the closed-loop system:
\bel{state*}\left\{\begin{aligned}
   d\h X(s) &=\big\{(A+\bar A+B\h\Psi_1)\h X+ (\ti A+B\h\Psi_2)\dbE[\h X]\big\}ds\\
   &\q+\big\{(C+\bar C+D\h\Psi_1)\h X+ (\ti C+D\h\Psi_2)\dbE[\h X] \big\}dW(s), \q s\in[0,T],\\
    \h X(0) &= \xi,
\end{aligned}\right.\ee
and $\h u(\cd)=\h\Psi_1(\cd) \h X(\cd)+\h\Psi_2(\cd)\dbE[\h X(\cd)]$ be the outcome of $(\h\Psi_1(\cd),\h\Psi_2(\cd))$.
We call $(\h\Psi_1(\cd),\h\Psi_2(\cd))$ a {\it na\"{\i}ve feedback strategy}, if
$$
\lim_{\|\D_N\|\to 0} \|\h u_{\D_N}(\cd)-\h u(\cd)\|=0,
$$
in which case, we call $\mathbf{\h\Psi}:=\h\Psi_1+\h \Psi_2\dbE$ a {\it na\"{\i}ve feedback operator}.
\end{definition}

\begin{theorem}\label{Thm:naive}
Let {\rm\ref{ass:A1}--\ref{ass:A2}}  hold. Define
$$
\h\Psi_1(s)=\Psi^*(s)+\bar\Psi^*(s),\q \h\Psi_2(s)=\ti\Psi^*(s),\q  s\in[0,T],
$$
where $(\Psi^*(\cd),\bar\Psi^*(\cd),\ti \Psi^*(\cd))$ is given by \rf{Prop:pre-solution-clsoed-main1}.
Then $(\h\Psi_1(\cd),\h\Psi_2(\cd))$ is a na\"{\i}ve feedback strategy and
$\mathbf{\h\Psi}=\h\Psi_1+\h \Psi_2\dbE$ is a na\"{\i}ve feedback operator.
\end{theorem}

The proof of Theorem \ref{Thm:naive} will be given in Section \ref{sec:nai}.

\subsection{Equilibrium Solutions}

The equilibrium solution  derives from a ``consistent planner" that optimizes the current problem with future controls fixed as constraints. Within this framework, at time $t$, the controller strategically engages with future decision-making selves by minimizing the cost functional over $[t, t+\e)$, while acknowledging control relinquishment beyond $t+\e$. The resulting solution is thus termed an equilibrium solution.

Motivated by \cite{Ekeland2010,Yong2012,Hu-Jin-Zhou2012,Bjork-Khapko-Murgoci2017},
we introduce the following definition of equilibrium solutions.

\begin{definition}\label{def:equilibrium}
We call $(\Psi_1^\dag(\cd),\Psi_2^\dag(\cd))\in L^\infty([0,T];\dbR^{m\times n})\times L^\infty([0,T];\dbR^{m\times n})$
an {\it equilibrium feedback strategy} if for any $t\in[0,T)$ and $u\in L^2_{\cF_t}(\Omega;\dbR^m)$, the following holds:
\bel{def:equilibrium1}
\liminf_{\e\to 0^+} {J(t,X^\dag(t);u^\e(\cd))-J(t,X^\dag (t);u^\dag(\cd))\over\e}\ges 0,
\ee
where
$$
u^\dag(s):=\Psi_1^\dag(s)X^\dag(s)+\Psi_2^\dag(s)\dbE[X^\dag(s)],\q s\in[0,T],
$$
with
$$\left\{\begin{aligned}
   dX^\dag(s) &=\big\{[A+ \bar A+B\Psi_1^\dag]X^\dag+[\ti A+B\Psi_2^\dag]\dbE[X^\dag] \big\}ds\\
   &\q+\big\{[C+ \bar C+D\Psi_1^\dag]X^\dag+[\ti C+D\Psi_2^\dag]\dbE[X^\dag] \big\}dW(s),\q s\in[0,T], \\
    X^\dag(0) &= \xi,
\end{aligned}\right.$$
and
$$
u^\e(s):=\left\{\begin{aligned}
&\Psi_1^\dag(s)X^\e(s)+\Psi_2^\dag(s)\dbE[X^\e(s)],\q &s\in[t+\e,T];\\
&\qq\qq\qq u,\q &s\in[t,t+\e),\end{aligned}\right.$$
with
$$\left\{\begin{aligned}
   dX^\e(s) &=\big\{A X^\e + \bar A\dbE_t[X^\e]+\ti A\dbE[X^\e]+B u \big\}ds\\
   &\q+\big\{C X^\e+ \bar C\dbE_t[X^\e]+\ti C\dbE[X^\e]+Du \big\}dW(s), \q s\in[t,t+\e);\\
    dX^\e(s) &=\big\{[A+ \bar A+B\Psi_1^\dag]X^\e+[\ti A+B\Psi_2^\dag]\dbE[X^\e] \big\}ds\\
   &\q+\big\{[C+ \bar C+D\Psi_1^\dag]X^\e+[\ti C+D\Psi_2^\dag]\dbE[X^\e] \big\}dW(s),\q s\in[t+\e,T],\\
    X^\e(t) &= X^\dag(t).
\end{aligned}\right.$$
If $(\Psi_1^\dag(\cd),\Psi_2^\dag(\cd))$ is an equilibrium feedback strategy, we call
$\mathbf{\Psi^\dag}=\Psi^\dag_1+ \Psi^\dag_2\dbE$ an {\it equilibrium feedback operator}.
\end{definition}

\begin{remark}
In Definition \ref{def:equilibrium}, the equilibrium strategy is defined via a variational method, introduced by \cite{Ekeland2010,Hu-Jin-Zhou2012,Bjork-Khapko-Murgoci2017}. An alternative definition of equilibrium strategies is provided by Yong \cite{Yong2012,Yong2014,Yong2017}, using a discretization approach. These two definitions are essentially equivalent. We will show that the equilibrium strategy constructed via Yong's multi-person differential game method \cite{Yong2012,Yong2014,Yong2017} satisfies Definition \ref{def:equilibrium}.
\end{remark}

To construct an equilibrium strategy, we introduce the following system of equilibrium Riccati equations:
\bel{ER-1}\left\{\begin{aligned}
  &\dot{\G}+\G\big [A+\bar A+B\Psi_1^\dag\big]+\big [A+\bar A+B\Psi_1^\dag\big]^\top \G\\
   &\q+\big [C+\bar C+D\Psi_1^\dag\big]^\top\G \big[C+\bar C+D\Psi_1^\dag\big]+Q+\Psi_1^{\dag\top} R\Psi_1^\dag=0, \\
   &\G(T)=G,
\end{aligned}\right.\ee
\bel{ER-2}\left\{\begin{aligned}
   &\dot{\bar\G}+\bar\G\big [A+\bar A+B\Psi_1^\dag\big]+\big [A+\bar A+B\Psi_1^\dag\big]^\top\bar \G\\
   &\q+\big [C+\bar C+D\Psi_1^\dag\big]^\top\G\big[C+\bar C+D\Psi_1^\dag\big]+Q+\bar Q+\Psi_1^{\dag\top} R\Psi_1^\dag=0, \\
   &\bar\G(T)=G+\bar G,
\end{aligned}\right.\ee
\bel{ER-3}\left\{\begin{aligned}
   &\dot{\ti\G}+\ti\G\big [A+\bar A+\ti A+B(\Psi_1^\dag+\Psi_2^\dag)\big]
   +\big[A+\bar A+\ti A+B(\Psi_1^\dag+\Psi_2^\dag)\big]^\top\ti \G\\
   &\q+\big [C+\bar C+\ti C+D(\Psi_1^\dag+\Psi_2^\dag)\big]^\top\G\big [C+\bar C+\ti C+D(\Psi_1^\dag+\Psi_2^\dag)\big]\\
   &\q +Q+\bar Q+\ti Q+(\Psi_1^\dag+\Psi_2^\dag)^\top R(\Psi_1^\dag+\Psi_2^\dag)=0, \\
   &\ti\G(T)=G+\bar G+\ti G,
\end{aligned}\right.\ee
where
\bel{ES1}
\Psi_1^\dag=\Psi^\dag+\bar\Psi^\dag,\q \Psi_2^\dag=\ti\Psi^\dag,
\ee
with
\begin{align}
\Psi^\dag&=-\big[R+D^\top \G D\big]^{-1}\big[B^\top \G+D^\top \G C],\nn\\
\bar\Psi^\dag&=-\big[R+D^\top \G D\big]^{-1}\big[B^\top(\bar\G- \G)+D^\top \G\bar C],\nn\\
\ti\Psi^\dag&=-\big[R+D^\top \G D\big]^{-1}\big[B^\top(\ti\G- \bar\G)+D^\top \G \ti C].
\label{ES}
\end{align}

\begin{theorem}\label{Thm:equiv}
Let {\rm\ref{ass:A1}--\ref{ass:A2}}  hold.
Then the system of equilibrium Riccati equations \rf{ER-1}--\rf{ER-3}
admits a unique solution $(\G(\cd),\bar\G(\cd),\ti \G(\cd))\in C([0,T];\dbS^n_+)^3$, and the function $(\Psi_1^\dag(\cd),\Psi^\dag_2(\cd))$
defined by \rf{ES1}--\rf{ES} is an equilibrium feedback strategy. Moreover, $\mathbf{\Psi^\dag}=\Psi^\dag_1+ \Psi^\dag_2\dbE$ is
an equilibrium feedback operator.
\end{theorem}

The well-posedness of \rf{ER-1}--\rf{ER-3} and the local optimality  \rf{def:equilibrium1}
will be proved in Section \ref{sec:equilib}
and Section \ref{sec:verifi}, respectively.

\section{Construction of  Na\"{\i}ve Solutions}\label{sec:nai}
Recall that
$\D_N=\{0=t_0<t_1<...<t_{N-1}<t_N=T\big\}$ is a partition of $[0,T]$ with
$\|\D_N\|=\sup_{1\les i\les N} |t_i-t_{i-1}|$.
At the time $t=t_0$, the state equation and  cost functional are given by
$$\left\{\begin{aligned}
   dX(s) &=\big\{AX+ \bar A\dbE[X]+\ti A\dbE[X]+Bu \big\}ds\\
   &\q+\big\{CX+ \bar C\dbE[X]+\ti C\dbE[X]+Du \big\}dW(s),\q s\in[t_0,T], \\
    X(t_0) &= \xi,
\end{aligned}\right.$$
and
\begin{align}
&J(t_0,\xi;u(\cd))
= \dbE\[\int_0^T \(\lan QX,X\ran +\lan\bar Q\dbE[X],\dbE[X]\ran
+\lan\ti Q\dbE[X],\dbE[X]\ran+\lan Ru,u\ran \) ds\nn\\
&\q +\lan GX(T),X(T)\ran +\lan\bar G\dbE[X(T)],\dbE[X(T)]\ran
+\lan\ti G\dbE[X(T)],\dbE[X(T)]\ran  \],\nn
\end{align}
respectively.
By Theorem \ref{Prop:pre-solution-clsoed},  the unique optimal control
$\h u_{\D_N}(\cd)$ on $[t_0,t_1)$ is given by
\begin{align*}
\h u_{\D_N}(s)=\Psi^*(s)\h X_{\D_N}(s)+\bar\Psi^*(s)\dbE[\h X_{\D_N}(s)]+\ti \Psi^*(s)\dbE[\h X_{\D_N}(s)],  \q s\in[t_0,t_1),
\end{align*}
with $\h X_{\D_N}(\cd)$ being the unique solution of the closed-loop system:
$$\left\{\begin{aligned}
   d\h X_{\D_N}(s) &=\big\{(A+B\Psi^*)\h X_{\D_N}+ (\bar A+B\bar\Psi^*)\dbE[\h X_{\D_N}]+ (\ti A+B\ti\Psi^*)\dbE[\h X_{\D_N}]\big\}ds\\
   &\q+\big\{(C+D\Psi^*)\h X_{\D_N}+ (\bar C+D\bar \Psi^*)\dbE[\h X_{\D_N}]\\
   &\q+ (\ti C+D\ti\Psi^*)\dbE[\h X_{\D_N}] \big\}dW(s), \q s\in[t_0,t_1),\\
    \h X_{\D_N}(t_0) &= \xi.
\end{aligned}\right.$$
Next, at the time $t=t_1$,
the state equation and the cost functional are given by
$$\left\{\begin{aligned}
   dX(s) &=\big\{AX+ \bar A\dbE_{t_1}[X]+\ti A\dbE[X]+Bu \big\}ds\\
   &\q+\big\{CX+ \bar C\dbE_{t_1}[X]+\ti C\dbE[X]+Du \big\}dW(s),\q s\in[t_1,T], \\
    X(t_1) &=\h X_{\D_N}(t_1) ,
\end{aligned}\right.$$
and
\begin{align*}
&J(t_1,\h X_{\D_N}(t_1) ;u(\cd))
= \dbE\[\int_{t_1}^T \(\lan QX,X\ran +\lan\bar Q\dbE_{t_1}[X],\dbE_{t_1}[X]\ran+\lan\ti Q\dbE[X],\dbE[X]\ran+\lan Ru,u\ran \) ds\\
&\q+\lan GX(T),X(T)\ran+\lan\bar G\dbE_{t_1}[X(T)],\dbE_{t_1}[X(T)]\ran
+\lan\ti G\dbE[X(T)],\dbE[X(T)]\ran  \],
\end{align*}
respectively.
By Theorem \ref{Prop:pre-solution-clsoed} again,  the unique optimal control
$\h u_{\D_N}(\cd)$ on $[t_1,t_2)$ is given by
\begin{align*}
\h u_{\D_N}(s)=\Psi^*(s)\h X_{\D_N}(s)+\bar\Psi^*(s)\dbE_{t_1}[\h X_{\D_N}(s)]+\ti \Psi^*(s)\dbE[\h X_{\D_N}(s)],  \q s\in[t_1,t_2),
\end{align*}
with $\h X_{\D_N}(\cd)$ being the unique solution of the closed-loop system:
$$\left\{\begin{aligned}
   d\h X_{\D_N}(s) &=\big\{(A+B\Psi^*)\h X_{\D_N}+ (\bar A+B\bar\Psi^*)\dbE_{t_1}[\h X_{\D_N}]+ (\ti A+B\ti\Psi^*)\dbE[\h X_{\D_N}]\big\}ds\\
   &\q+\big\{(C+D\Psi^*)\h X_{\D_N}+ (\bar C+D\bar \Psi^*)\dbE_{t_1}[\h X_{\D_N}]\\
   &\qq+ (\ti C+D\ti\Psi^*)\dbE[\h X_{\D_N}] \big\}dW(s), \q s\in[t_1,t_2),\\
    \h X_{\D_N}(t_1) &=\h X_{\D_N}(t_1).
\end{aligned}\right.$$
By induction, the na\"{\i}ve control $\h u_{\D_N}(\cd)$   associated with $\D_N$ is given by
\begin{align}
\h u_{\D_N}(s)&=\Psi^*(s)\h X_{\D_N}(s)+\bar\Psi^*(s)\dbE_{t_i}[\h X_{\D_N}(s)]+\ti \Psi^*(s)\dbE[\h X_{\D_N}(s)],\nn  \\
&\qq\qq\qq\qq\qq s\in[t_i,t_{i+1}),\q i=0,1,...,N-2,\nn\\
\h u_{\D_N}(s)&=\Psi^*(s)\h X_{\D_N}(s)+\bar\Psi^*(s)\dbE_{t_{N-1}}[\h X_{\D_N}(s)]+\ti \Psi^*(s)\dbE[\h X_{\D_N}(s)],\nn\\
&\qq\qq\qq\qq\qq s\in[t_{N-1},t_N],
\end{align}
with $\h X_{\D_N}(\cd)$ being the unique solution of the closed-loop system:
\bel{nai-XN}\left\{\begin{aligned}
   d\h X_{\D_N}(s) &=\big\{(A+B\Psi^*)\h X_{\D_N}+ (\bar A+B\bar\Psi^*)\dbE_{t_i}[\h X_{\D_N}]+ (\ti A+B\ti\Psi^*)\dbE[\h X_{\D_N}]\big\}ds\\
   &\q+\big\{(C+D\Psi^*)\h X_{\D_N}+ (\bar C+D\bar \Psi^*)\dbE_{t_i}[\h X_{\D_N}]\\
   &\qq+ (\ti C+D\ti\Psi^*)\dbE[\h X_{\D_N}] \big\}dW(s), \q s\in[t_i,t_{i+1}),\q i=0,1,...,N-2;\\
    d\h X_{\D_N}(s) &=\big\{(A+B\Psi^*)\h X_{\D_N}+ (\bar A+B\bar\Psi^*)\dbE_{t_{N-1}}[\h X_{\D_N}]+ (\ti A+B\ti\Psi^*)\dbE[\h X_{\D_N}]\big\}ds\\
   &\q+\big\{(C+D\Psi^*)\h X_{\D_N}+ (\bar C+D\bar \Psi^*)\dbE_{t_{N-1}}[\h X_{\D_N}]\\
   &\qq+ (\ti C+D\ti\Psi^*)\dbE[\h X_{\D_N}] \big\}dW(s), \q s\in[t_{N-1},t_N],\\
   \h X_{\D_N}(0) &= \xi.
\end{aligned}\right.\ee
Recall from \rf{state*} that
\bel{state-h}\left\{\begin{aligned}
   d\h X(s) &=\big\{(A+\bar A+B\Psi^*+B\bar\Psi^*)\h X+ (\ti A+B\ti\Psi^*)\dbE[\h X]\big\}ds\\
   &\q+\big\{(C+\bar C+D\Psi^*+D\bar \Psi^*)\h X+ (\ti C+D\ti\Psi^*)\dbE[\h X] \big\}dW(s), \q s\in[0,T],\\
    \h X(0) &= \xi.
\end{aligned}\right.\ee
By applying the standard estimates of SDEs to \rf{nai-XN}--\rf{state-h}, we have
\bel{state-h-est}
\dbE\[\sup_{s\in[0,T]}|\h X_{\D_N}(s)-\h X(s)|^2\]\les K \sum_{i=0}^{N-1}\dbE\[\int_{t_i}^{t_{i+1}}\big|\h X(s)-\dbE_{t_i}[\h X(s)]\big|^2ds\].
\ee
Note that for $s\in[t_i,t_{i+1}]$, we have
\begin{align*}
&\h X(s)-\dbE_{t_i}[\h X(s)]=\int_{t_i}^s(A+\bar A+B\Psi^*+B\bar\Psi^*)(\h X-\dbE_{t_i}[\h X])dr\\
&\qq+\int_{t_i}^s\big\{(C+\bar C+D\Psi^*+D\bar \Psi^*)\h X+ (\ti C+D\ti\Psi^*)\dbE[\h X] \big\}dW(r).
\end{align*}
Thus,
\begin{align*}
&\dbE\big[|\h X(s)-\dbE_{t_i}[\h X(s)]|^2\big]\les K \dbE\[\sup_{r\in[0,T]}|\h X(r)|^2\]|s-t_i|\les K\|\D_N\|,\q s\in[t_i,t_{i+1}].
\end{align*}
Substituting the above into \rf{state-h-est} yields that
\begin{align}
\lim_{\|\D_N\| \to0}\dbE\[\sup_{s\in[0,T]}|\h X_{\D_N}(s)-\h X(s)|^2\]\les\lim_{\|\D_N\| \to0}KT\|\D_N\|=0.
\end{align}
Moreover,
\begin{align}
&\dbE\[\int_0^T |\h u_{\D_N}(s)-\h u(s)|^2ds\]\les K\dbE\[\int_0^T |\h X_{\D_N}(s)-\h X(s)|^2ds\]\nn\\
&\qq +K \sum_{i=0}^{N-1}\dbE\[\int_{t_i}^{t_{i+1}}\big|\h X(s)-\dbE_{t_i}[\h X_{\D_N}(s)]\big|^2ds\]\nn\\
&\q\les K\dbE\[\int_0^T |\h X_{\D_N}(s)-\h X(s)|^2ds\]+K \sum_{i=0}^{N-1}\dbE\[\int_{t_i}^{t_{i+1}}\big|\h X(s)-\dbE_{t_i}[ \h X(s)]\big|^2ds\]\nn\\
&\qq+K \sum_{i=0}^{N-1}\dbE\[\int_{t_i}^{t_{i+1}}\big|\dbE_{t_i}[ \h X(s)]-\dbE_{t_i}[ \h X_{\D_N}(s)]\big|^2ds\]\nn\\
&\q\les K\dbE\[\int_0^T |\h X_{\D_N}(s)-\h X(s)|^2ds\]+K \sum_{i=0}^{N-1}\dbE\[\int_{t_i}^{t_{i+1}}\big|\h X(s)-\dbE_{t_i}[ \h X(s)]\big|^2ds\]\nn\\
&\q\to0,\q\hbox{as $\|\D_N\| \to0$}.
\end{align}
This completes the proof of Theorem \ref{Thm:naive}.

\section{Construction of  Equilibrium Solutions: A Multi-Person Differential Game Approach}\label{sec:equilib}

Let $\D_N=\{0=t_0<t_1<...<t_{N-1}<t_N=T\big\}$ be a partition of $[0,T]$ with
$\|\D_N\|=\sup_{1\les i\les N} |t_i-t_{i-1}|$.
For simplicity, we usually write $\D_N$ as $\D$ if there is no confusion.

\ms

\textbf{Problem on $[t_{N-1},T]$.}
At the time $t=t_{N-1}$, the state equation and  cost functional are given by
$$\left\{\begin{aligned}
   &dX_{N-1}(s) =\big\{AX_{N-1}+ \bar A\dbE_{t_{N-1}}[X_{N-1}]+\ti A\dbE[X_{N-1}]+Bu_{N-1} \big\}ds\\
   &\qq\qq\qq+\big\{CX_{N-1}+ \bar C\dbE_{t_{N-1}}[X_{N-1}]+\ti C\dbE[X_{N-1}]+Du_{N-1} \big\}dW(s),\q s\in[t_{N-1},T], \\
    &X_{N-1}(t_{N-1}) = \xi,
\end{aligned}\right.$$
and
\begin{align}
&J_{N-1}(t_{N-1},\xi;u(\cd))
= \dbE\[\int_{t_{N-1}}^T \(\lan QX_{N-1},X_{N-1}\ran +\lan\bar Q\dbE_{t_{N-1}}[X_{N-1}],\dbE_{t_{N-1}}[X_{N-1}]\ran\nn\\
&\q+\lan\ti Q\dbE[X_{N-1}],\dbE[X_{N-1}]\ran+\lan Ru_{N-1},u_{N-1}\ran \) ds+\lan GX_{N-1}(T),X_{N-1}(T)\ran\nn\\
&\q  +\lan\bar G\dbE_{t_{N-1}}[X_{N-1}(T)],\dbE_{t_{N-1}}[X_{N-1}(T)]\ran
+\lan\ti G\dbE[X_{N-1}(T)],\dbE[X_{N-1}(T)]\ran  \],\nn
\end{align}
respectively.
By Theorem \ref{Prop:pre-solution-clsoed},  the unique optimal control
$ u^\dag_{N-1}(\cd)$ on $[t_{N-1},t_N]$ is given by
\begin{align*}
u^\dag_{N-1} (s)&=\Psi^\dag_{N-1}(s)X^\dag_{N-1}(s)+\bar\Psi^\dag_{N-1}(s)\dbE_{t_{N-1}}[X^\dag_{N-1}(s)]\\
&\q+\ti \Psi^\dag_{N-1}(s)\dbE[X^\dag_{N-1}(s)],  \q s\in[t_{N-1},t_N],
\end{align*}
where
\bel{state-dag-(N-1)}\left\{\begin{aligned}
  & dX^\dag_{N-1}(s) =\big\{(A+B\Psi^\dag_{N-1})X^\dag_{N-1}+ (\bar A+B\bar\Psi^\dag_{N-1})\dbE_{t_{N-1}}[X^\dag_{N-1}]\\
   &\qq\qq\q~+ (\ti A+B\ti\Psi^\dag_{N-1})\dbE[X^\dag_{N-1}]\big\}ds\\
   &\qq\qq\q~+\big\{(C+D\Psi^\dag_{N-1})X^\dag_{N-1}+ (\bar C+D\bar\Psi^\dag_{N-1})\dbE_{t_{N-1}}[X^\dag_{N-1}]\\
   &\qq\qq\q~+ (\ti C+D\ti\Psi^\dag_{N-1})\dbE[X^\dag_{N-1}]\big\}dW(s), \q s\in[t_{N-1},t_N],\\
    &X^\dag_{N-1}(t_{N-1}) = \xi,
\end{aligned}\right.\ee
and
\begin{align*}
\Psi_{N-1}^\dag&=-\big[R+D^\top P_{N-1}D\big]^{-1}\big[B^\top P_{N-1}+D^\top P_{N-1}C\big],\\
\bar\Psi_{N-1}^\dag&=-\big[R+D^\top P_{N-1}D\big]^{-1}\big[B^\top(\Pi_{N-1}- P_{N-1})+D^\top P_{N-1}\bar C\big],\\
\ti\Psi_{N-1}^\dag&=-\big[R+D^\top P_{N-1}D\big]^{-1}\big[B^\top(\Phi_{N-1}- \Pi_{N-1})+D^\top P_{N-1}\ti C\big],
\end{align*}
with the Riccati equations over $[t_{N-1},t_N]$:
$$\left\{\begin{aligned}
   &\dot{P}_{N-1}+P_{N-1}A+A^\top P_{N-1}+C^\top P_{N-1}C+Q\\
   &\q-(P_{N-1}B+C^\top P_{N-1}D)(R+D^\top P_{N-1}D)^{-1}(B^\top P_{N-1}+D^\top P_{N-1}C)=0, \\
   &P_{N-1}(T)=G,
\end{aligned}\right.$$
$$\left\{\begin{aligned}
   &\dot{\Pi}_{N-1}+\Pi_{N-1}(A+\bar A)+(A+\bar A)^\top \Pi_{N-1}+(C+\bar C)^\top P_{N-1}(C+\bar C)+Q+\bar Q\\
   &\q-[\Pi_{N-1} B+(C+\bar C)^\top P_{N-1}D](R+D^\top P_{N-1}D)^{-1}[B^\top \Pi_{N-1}+D^\top P_{N-1}(C+\bar C)]=0,\\
   &\Pi_{N-1}(T)=G+\bar G,
\end{aligned}\right.$$
$$
\left\{\begin{aligned}
  &\dot{\Phi}_{N-1}+\Phi_{N-1}(A+\bar A+\ti A)+(A+\bar A+\ti A)^\top \Phi_{N-1}\\
   &\qq+(C+\bar C+\ti C)^\top P_{N-1}(C+\bar C+\ti C)+Q+\bar Q+\ti Q\\
   &\qq-[\Phi_{N-1} B+(C+\bar C+\ti C)^\top P_{N-1}D](R+D^\top P_{N-1}D)^{-1}\\
   &\qq\times[B^\top \Phi_{N-1}+D^\top P_{N-1}(C+\bar C+\ti C)]=0,\\
   &\Phi_{N-1}(T)=G+\bar G+\ti G.
\end{aligned}\right.
$$

\ms
\textbf{Problem on $[t_{N-2},T]$.} On $[t_{N-2},T]$, we have the state equation
\bel{state-(N-2)}\left\{\begin{aligned}
   &dX_{N-2}(s) =\big\{AX_{N-2}+ \bar A\dbE_{t_{N-2}}[X_{N-2}]+\ti A\dbE[X_{N-2}]+Bu_{N-2}\big\}ds\\
   &\qq\qq\qq+\big\{CX_{N-2}+ \bar C\dbE_{t_{N-2}}[X_{N-2}]+\ti C\dbE[X_{N-2}]+Du_{N-2} \big\}dW(s),\\
   &\qq\qq\qq\qq s\in[t_{N-2},t_{N-1}), \\
    &X_{N-2}(t_{N-2}) = \xi,
\end{aligned}\right.\ee
and the cost functional
\begin{align}\label{cost-(N-2)}
&J_{N-2}(t_{N-2},\xi;u(\cd))
= \dbE\Big\{\int_{t_{N-2}}^{t_{N-1}} \[\lan Q X_{N-2},X_{N-2}\ran +\lan\bar Q\dbE_{t_{N-2}}[X_{N-2}],\dbE_{t_{N-2}}[X_{N-2}]\ran\nn\\
&\qq +\lan\ti Q\dbE[X_{N-2}],\dbE[X_{N-2}]\ran+\lan Ru_{N-2},u_{N-2}\ran \] ds
+\int_{t_{N-1}}^{T} \[\lan Q X^\dag_{N-1},X_{N-1}^\dag\ran \nn\\
&\qq+\lan\bar Q\dbE_{t_{N-2}}[X_{N-1}^\dag],\dbE_{t_{N-2}}[X_{N-1}^\dag]\ran
+\lan\ti Q\dbE[X_{N-1}^\dag],\dbE[X_{N-1}^\dag]\ran+\lan Ru^\dag_{N-1},u^\dag_{N-1}\ran \] ds\nn\\
&\qq+\lan GX^\dag_{N-1}(T),X^\dag_{N-1}(T)\ran+\lan\bar G\dbE_{t_{N-2}}[X^\dag_{N-1}(T)],\dbE_{t_{N-2}}[X^\dag_{N-1}(T)]\ran\nn\\
&\qq+\lan\ti G\dbE[X^\dag_{N-1}(T)],\dbE[X^\dag_{N-1}(T)]\ran  \Big\},
\end{align}
in which $X^\dag_{N-1}(\cd)$ is the unique solution of \rf{state-dag-(N-1)} with the initial state $X^\dag_{N-1}(t_{N-1})=X_{N-2}(t_{N-1})$.
Note that
\begin{align}
&\dbE\big[\lan Ru^\dag_{N-1},u^\dag_{N-1}\ran\big]
=\dbE\[\Big\lan R\big(\Psi^\dag_{N-1}X^\dag_{N-1}+\bar\Psi^\dag_{N-1}\dbE_{t_{N-1}}[X^\dag_{N-1}]+\ti \Psi^\dag_{N-1}\dbE[X^\dag_{N-1}]\big),\nn\\
&\qq\q\big(\Psi^\dag_{N-1}X^\dag_{N-1}+\bar\Psi^\dag_{N-1}\dbE_{t_{N-1}}[X^\dag_{N-1}]+\ti \Psi^\dag_{N-1}\dbE[X^\dag_{N-1}]\big)\Big\ran\]\nn\\
&=\dbE\[\big\lan\big (\Psi_{N-1}^{\dag\top} R\Psi_{N-1}^\dag\big)X^\dag_{N-1},\,X^\dag_{N-1}\big\ran
+\big\lan\big (\Psi_{N-1}^{\dag\top} R\bar\Psi_{N-1}^\dag+\bar\Psi_{N-1}^{\dag\top} R\Psi_{N-1}^\dag\nn\\
&\qq +\bar\Psi_{N-1}^{\dag\top} R\bar\Psi_{N-1}^\dag\big)\dbE_{t_{N-1}}[X_{N-1}^\dag],\,\dbE_{t_{N-1}}[X_{N-1}^\dag]\big\ran
+\big\lan\big (\Psi_{N-1}^{\dag\top} R\ti\Psi_{N-1}^\dag+\ti\Psi_{N-1}^{\dag\top} R\Psi_{N-1}^\dag\nn\\
&\qq
+\ti\Psi_{N-1}^{\dag\top} R\bar\Psi_{N-1}^\dag
+\bar\Psi_{N-1}^{\dag\top} R\ti\Psi_{N-1}^\dag+\ti\Psi_{N-1}^{\dag\top} R\ti\Psi_{N-1}^\dag\big)\dbE[X_{N-1}^\dag],\,\dbE[X_{N-1}^\dag]\big\ran\].\nn
\end{align}
Then
\begin{align}
\dbI_{N-1}&:=\dbE\Big\{\int_{t_{N-1}}^{T} \[\lan Q X^\dag_{N-1},X_{N-1}^\dag\ran+\lan\bar Q\dbE_{t_{N-2}}[X_{N-1}^\dag],\dbE_{t_{N-2}}[X_{N-1}^\dag]\ran
 \nn\\
&\qq+\lan\ti Q\dbE[X_{N-1}^\dag],\dbE[X_{N-1}^\dag]\ran+\lan Ru^\dag_{N-1},u^\dag_{N-1}\ran \] ds+\lan GX^\dag_{N-1}(T),X^\dag_{N-1}(T)\ran\nn\\
&\qq+\lan\bar G\dbE_{t_{N-2}}[X^\dag_{N-1}(T)],\dbE_{t_{N-2}}[X^\dag_{N-1}(T)]\ran+\lan\ti G\dbE[X^\dag_{N-1}(T)],\dbE[X^\dag_{N-1}(T)]\ran  \Big\}\nn\\
&=\dbE\Big\{\int_{t_{N-1}}^{T}\[\big\lan\big (Q+\Psi_{N-1}^{\dag\top} R\Psi_{N-1}^\dag\big)X^\dag_{N-1},\,X^\dag_{N-1}\big\ran
\nn\\
&\qq +\big\lan\big (\Psi_{N-1}^{\dag\top} R\bar\Psi_{N-1}^\dag+\bar\Psi_{N-1}^{\dag\top} R\Psi_{N-1}^\dag+\bar\Psi_{N-1}^{\dag\top} R\bar\Psi_{N-1}^\dag\big)\dbE_{t_{N-1}}[X_{N-1}^\dag],\,\dbE_{t_{N-1}}[X_{N-1}^\dag]\big\ran
\nn\\
&\qq +\big\lan\bar Q\dbE_{t_{N-2}}[X_{N-1}^\dag],\,\dbE_{t_{N-2}}[X_{N-1}^\dag]\big\ran
+\big\lan\big (\ti Q+\Psi_{N-1}^{\dag\top} R\ti\Psi_{N-1}^\dag+\ti\Psi_{N-1}^{\dag\top} R\Psi_{N-1}^\dag\nn\\
&\qq+\ti\Psi_{N-1}^{\dag\top} R\bar\Psi_{N-1}^\dag
+\bar\Psi_{N-1}^{\dag\top} R\ti\Psi_{N-1}^\dag+\ti\Psi_{N-1}^{\dag\top} R\ti\Psi_{N-1}^\dag\big)\dbE[X_{N-1}^\dag],\,\dbE[X_{N-1}^\dag]\big\ran\]ds\nn\\
&\qq+\lan GX^\dag_{N-1}(T),X^\dag_{N-1}(T)\ran+\lan\bar G\dbE_{t_{N-2}}[X^\dag_{N-1}(T)],\dbE_{t_{N-2}}[X^\dag_{N-1}(T)]\ran
  \nn\\
&\qq+\lan\ti G\dbE[X^\dag_{N-1}(T)],\dbE[X^\dag_{N-1}(T)]\ran  \Big\}.\nn
\end{align}
Then by Lemma \ref{lem:rep}, we have
\begin{align}
\dbI_{N-1}&=\dbE\[\big\lan\G_{N-2}(t_{N-1})X_{N-2}(t_{N-1}),\,X_{N-2}(t_{N-1})\big\ran\nn\\
&\q+\big\lan [\bar\G_{N-2}(t_{N-1})-\G_{N-2}(t_{N-1})]\dbE_{t_{N-2}}[X_{N-2}(t_{N-1})],\,\dbE_{t_{N-2}}[X_{N-2}(t_{N-1})]\big\ran\nn\\
&\q+\big\lan [\ti\G_{N-2}(t_{N-1})-\bar  \G_{N-2}(t_{N-1})]\dbE[X_{N-2}(t_{N-1})],\,\dbE[X_{N-2}(t_{N-1})]\big\ran\],\label{repre-(N-1)}
\end{align}
with the Lyapunov equations over $[t_{N-1},T]$:
$$\left\{\begin{aligned}
   &\dot{\check\G}_{N-2}+\check\G_{N-2} (A+B\Psi_{N-1}^\dag)+(A+B\Psi_{N-1}^\dag)^\top \check\G_{N-2}\\
   &\q+ (C+D\Psi_{N-1}^\dag)^\top \check\G_{N-2}  (C+D\Psi_{N-1}^\dag)+Q+\Psi_{N-1}^{^\dag\top} R\Psi_{N-1}^\dag=0, \\
   &\check\G_{N-2}(T)=G,
\end{aligned}\right.$$
$$\left\{\begin{aligned}
  &\dot{\G}_{N-2}+\G_{N-2} [A+\bar A+B(\Psi_{N-1}^\dag+\bar\Psi_{N-1}^\dag)]\\
   &\q+ [A+\bar A+B(\Psi_{N-1}^\dag+\bar\Psi_{N-1}^\dag)]^\top \G_{N-2}\\
   &\q+ [C+\bar C+D(\Psi_{N-1}^\dag+\bar\Psi_{N-1}^\dag)]^\top\check\G_{N-2} [C+\bar C+D(\Psi_{N-1}^\dag+\bar\Psi_{N-1}^\dag)]\\
   &\q+Q+(\Psi_{N-1}^\dag+\bar\Psi_{N-1}^\dag)^\top R(\Psi_{N-1}^\dag+\bar\Psi_{N-1}^\dag)=0, \\
   &\G_{N-2}(T)=G,
\end{aligned}\right.$$
$$\left\{\begin{aligned}
   &\dot{\bar\G}_{N-2}+\bar\G_{N-2} [A+\bar A+B(\Psi_{N-1}^\dag+\bar\Psi_{N-1}^\dag)]\\
   &\q+ [A+\bar A+B(\Psi_{N-1}^\dag+\bar\Psi_{N-1}^\dag)]^\top\bar \G_{N-2}\\
   &\q+ [C+\bar C+D(\Psi_{N-1}^\dag+\bar\Psi_{N-1}^\dag)]^\top\check\G_{N-2} [C+\bar C+D(\Psi_{N-1}^\dag+\bar\Psi_{N-1}^\dag)]\\
   &\q+Q+\bar Q+(\Psi_{N-1}^\dag+\bar\Psi_{N-1}^\dag)^\top R(\Psi_{N-1}^\dag+\bar\Psi_{N-1}^\dag)=0, \\
   &\bar\G_{N-2}(T)=G+\bar G,
\end{aligned}\right.$$
$$\left\{\begin{aligned}
   &\dot{\ti\G}_{N-2}+\ti\G_{N-2} [A+\bar A+\ti A+B(\Psi_{N-1}^\dag+\bar\Psi_{N-1}^\dag+\ti\Psi_{N-1}^\dag)]\\
   &\q+[A+\bar A+\ti A+B(\Psi_{N-1}^\dag+\bar\Psi_{N-1}^\dag+\ti\Psi_{N-1}^\dag)]^\top\ti \G_{N-2}\\
   &\q+ [C+\bar C+\ti C+D(\Psi_{N-1}^\dag+\bar\Psi_{N-1}^\dag+\ti\Psi_{N-1}^\dag)]^\top\check\G_{N-2}\\
   &\q\times  [C+\bar C+\ti C+D(\Psi_{N-1}^\dag+\bar\Psi_{N-1}^\dag+\ti\Psi_{N-1}^\dag)]+Q+\bar Q+\ti Q\\
   &\q+(\Psi_{N-1}^\dag+\bar\Psi_{N-1}^\dag+\ti\Psi_{N-1}^\dag)^\top R(\Psi_{N-1}^\dag+\bar\Psi_{N-1}^\dag+\ti\Psi_{N-1}^\dag)=0, \\
   &\ti\G_{N-2}(T)=G+\bar G+\ti G.
\end{aligned}\right.$$
Substituting \rf{repre-(N-1)} into \rf{cost-(N-2)} yields that
\begin{align}\label{cost-(N-2)-1}
&J_{N-2}(t_{N-2},\xi;u(\cd))
= \dbE\Big\{\int_{t_{N-2}}^{t_{N-1}} \[\lan Q X_{N-2},X_{N-2}\ran \nn\\
&\qq+\lan\bar Q\dbE_{t_{N-2}}[X_{N-2}],\dbE_{t_{N-2}}[X_{N-2}]\ran +\lan\ti Q\dbE[X_{N-2}],\dbE[X_{N-2}]\ran\nn\\
&\qq+\lan Ru_{N-2},u_{N-2}\ran \] ds+\big\lan\G_{N-2}(t_{N-1})X_{N-2}(t_{N-1}),\,X_{N-2}(t_{N-1})\big\ran\nn\\
&\qq+\big\lan [\bar\G_{N-2}(t_{N-1})-\G_{N-2}(t_{N-1})]\dbE_{t_{N-2}}[X_{N-2}(t_{N-1})],\,\dbE_{t_{N-2}}[X_{N-2}(t_{N-1})]\big\ran\nn\\
&\qq+\big\lan [\ti\G_{N-2}(t_{N-1})-\bar  \G_{N-2}(t_{N-1})]\dbE[X_{N-2}(t_{N-1})],\,\dbE[X_{N-2}(t_{N-1})]\big\ran  \Big\}.
\end{align}
Then we end with a standard mean-field SLQ optimal control problem over $[t_{N-2},t_{N-1}]$ with the state equation \rf{state-(N-2)}
and  cost functional \rf{cost-(N-2)-1}.
By Theorem \ref{Prop:pre-solution-clsoed} again,  the unique optimal control
$ u^\dag_{N-2}(\cd)$ on $[t_{N-2},t_{N-1})$ is given by
\begin{align*}
u^\dag_{N-2} (s)&=\Psi^\dag_{N-2}(s)X^\dag_{N-2}(s)+\bar\Psi^\dag_{N-2}(s)\dbE_{t_{N-2}}[X^\dag_{N-2}(s)]\\
&\q+\ti \Psi^\dag_{N-2}(s)\dbE[X^\dag_{N-2}(s)],  \q s\in[t_{N-2},t_{N-1}),
\end{align*}
where
$$\left\{\begin{aligned}
  & dX^\dag_{N-2}(s) =\big\{(A+B\Psi^\dag_{N-2})X^\dag_{N-2}+ (\bar A+B\bar\Psi^\dag_{N-2})\dbE_{t_{N-2}}[X^\dag_{N-2}]\\
   &\qq\qq\q~+ (\ti A+B\ti\Psi^\dag_{N-2})\dbE[X^\dag_{N-2}]\big\}ds\\
   &\qq\qq\q~+\big\{(C+D\Psi^\dag_{N-2})X^\dag_{N-2}+ (\bar C+D\bar\Psi^\dag_{N-2})\dbE_{t_{N-2}}[X^\dag_{N-2}]\\
   &\qq\qq\q~+ (\ti C+D\ti\Psi^\dag_{N-2})\dbE[X^\dag_{N-2}]\big\}dW(s), \q s\in[t_{N-2},t_{N-1}),\\
    &X^\dag_{N-2}(t_{N-2}) = \xi,
\end{aligned}\right.$$
and
\begin{align*}
\Psi_{N-2}^\dag&=-\big[R+D^\top P_{N-2}D\big]^{-1}\big[B^\top P_{N-2}+D^\top P_{N-2}C\big],\\
\bar\Psi_{N-2}^\dag&=-\big[R+D^\top P_{N-2}D\big]^{-1}\big[B^\top(\Pi_{N-2}- P_{N-2})+D^\top P_{N-2}\bar C\big],\\
\ti\Psi_{N-2}^\dag&=-\big[R+D^\top P_{N-2}D\big]^{-1}\big[B^\top(\Phi_{N-2}- \Pi_{N-2})+D^\top P_{N-2}\ti C\big],
\end{align*}
with the Riccati equations over $[t_{N-2},t_{N-1}]$:
\begin{equation*}\left\{\begin{aligned}
   &\dot{P}_{N-2}+P_{N-2}A+A^\top P_{N-2}+C^\top P_{N-2}C+Q\\
   &\q-(P_{N-2}B+C^\top P_{N-2}D)(R+D^\top P_{N-2}D)^{-1}(B^\top P_{N-2}+D^\top P_{N-2}C)=0, \\
   &P_{N-2}(t_{N-1})=\G_{N-2}(t_{N-1}),
\end{aligned}\right.\end{equation*}
\begin{equation*}\left\{\begin{aligned}
   &\dot{\Pi}_{N-2}+\Pi_{N-2}(A+\bar A)+(A+\bar A)^\top \Pi_{N-2}+(C+\bar C)^\top P_{N-2}(C+\bar C)+Q+\bar Q\\
   &\q-[\Pi_{N-2} B+(C+\bar C)^\top P_{N-2}D](R+D^\top P_{N-2}D)^{-1}[B^\top \Pi_{N-2}+D^\top P_{N-2}(C+\bar C)]=0,\\
   &\Pi_{N-2}(t_{N-1})=\bar\G_{N-2}(t_{N-1}),
\end{aligned}\right.\end{equation*}
\begin{equation*}\left\{\begin{aligned}
  &\dot{\Phi}_{N-2}+\Phi_{N-2}(A+\bar A+\ti A)+(A+\bar A+\ti A)^\top \Phi_{N-2}\\
   &\qq+(C+\bar C+\ti C)^\top P_{N-2}(C+\bar C+\ti C)+Q+\bar Q+\ti Q\\
   &\qq-[\Phi_{N-2} B+(C+\bar C+\ti C)^\top P_{N-2}D](R+D^\top P_{N-2}D)^{-1}\\
   &\qq\times[B^\top \Phi_{N-2}+D^\top P_{N-2}(C+\bar C+\ti C)]=0,\\
   &\Phi_{N-2}(t_{N-1})=\ti\G_{N-2}(t_{N-1}).
\end{aligned}\right.\end{equation*}
Define
\begin{equation*}
\f_\D(s)=\left\{\begin{aligned}
  &\f_{N-1}(s),\q s\in(t_{N-1},T];\\
   &\f_{N-2}(s),\q s\in[t_{N-2},t_{N-1}],
\end{aligned}\right.
\end{equation*}
for $\a(\cd)=P(\cd),\Pi(\cd),\Phi(\cd),X^\dag(\cd),u^\dag(\cd),\Psi^\dag(\cd),\bar\Psi^\dag(\cd),\ti\Psi^\dag(\cd)$.
 Moreover, we define the function $\rho_\D(\cd)$ by
\begin{equation*}
\rho_\D(s)=\left\{\begin{aligned}
  &t_{N-1},\q s\in[t_{N-1},T];\\
   &\q t_i,\q~\, s\in[t_i,t_{i+1}),\q i=0,1,...,N-2.
\end{aligned}\right.
\end{equation*}
Then on $[t_{N-2},T]$ we have
\bel{state-D-(N-2)}\left\{\begin{aligned}
  & dX^\dag_\D(s) =\big\{(A+B\Psi^\dag_\D)X^\dag_\D+ (\bar A+B\bar\Psi^\dag_\D)\dbE_{\rho_\D(s)}[X^\dag_\D]+ (\ti A+B\ti\Psi^\dag_\D)\dbE[X^\dag_\D]\big\}ds\\
   &\qq\qq~+\big\{(C+D\Psi^\dag_\D)X^\dag_\D+ (\bar C+D\bar\Psi^\dag_\D)\dbE_{\rho_\D(s)}[X^\dag_\D]\\
   &\qq\qq~+ (\ti C+D\ti\Psi^\dag_\D)\dbE[X^\dag_\D]\big\}dW(s), \q s\in[t_{N-2},T],\\
    &X^\dag_\D(t_{N-2}) = \xi,
\end{aligned}\right.\ee
with the control process:
\begin{align*}
u^\dag_\D (s)&=\Psi^\dag_\D(s)X^\dag_\D(s)+\bar\Psi^\dag_\D(s)\dbE_{\rho_\D(s)}[X^\dag_\D(s)]+\ti \Psi^\dag_\D(s)\dbE[X^\dag_\D(s)],  \q s\in[t_{N-2},T].
\end{align*}

\ms
\textbf{Problem on $[t_{N-3},T]$.}
On $[t_{N-3},t_N]$, we have the state equation
\bel{state-(N-3)}\left\{\begin{aligned}
   &dX_{N-3}(s) =\big\{AX_{N-3}+ \bar A\dbE_{t_{N-3}}[X_{N-3}]+\ti A\dbE[X_{N-3}]+Bu_{N-3} \big\}ds\\
   &\qq\qq\qq+\big\{CX_{N-3}+ \bar C\dbE_{t_{N-3}}[X_{N-3}]+\ti C\dbE[X_{N-3}]+Du_{N-3}\big\}dW(s),\\
   &\qq\qq\qq\qq s\in[t_{N-3},t_{N-2}), \\
    &X_{N-3}(t_{N-3}) = \xi,
\end{aligned}\right.\ee
and the cost functional
\begin{align}
&J_{N-3}(t_{N-3},\xi;u(\cd))
= \dbE\Big\{\int_{t_{N-3}}^{t_{N-2}} \[\lan Q X_{N-3},X_{N-3}\ran+\lan\bar Q\dbE_{t_{N-3}}[X_{N-3}],\dbE_{t_{N-3}}[X_{N-3}]\ran \nn\\
&\qq +\lan\ti Q\dbE[X_{N-3}],\dbE[X_{N-3}]\ran+\lan Ru_{N-3},u_{N-3}\ran \] ds
+\int_{t_{N-2}}^{T} \[\lan Q X^\dag_\D,X^\dag_\D\ran\nn\\
&\qq +\lan\bar Q\dbE_{t_{N-3}}[X^\dag_\D],\dbE_{t_{N-3}}[X^\dag_\D]\ran
+\lan\ti Q\dbE[X^\dag_\D],\dbE[X^\dag_\D]\ran+\lan Ru^\dag_\D,u^\dag_\D\ran \] ds\nn\\
&\qq+\lan GX_\D(T),X_\D(T)\ran+\lan\bar G\dbE_{t_{N-3}}[X_\D(T)],\dbE_{t_{N-3}}[X_\D(T)]\ran\nn\\
&\qq+\lan\ti G\dbE[X_{\D}(T)],\dbE[X_{\D}(T)]\ran  \Big\}\nn\\
&=: \dbE\Big\{\int_{t_{N-3}}^{t_{N-2}} \[\lan Q X_{N-3},X_{N-3}\ran  +\lan\bar Q\dbE_{t_{N-3}}[X_{N-3}],\dbE_{t_{N-3}}[X_{N-3}]\ran\nn\\
&\qq+\lan\ti Q\dbE[X_{N-3}],\dbE[X_{N-3}]\ran+\lan Ru_{N-3},u_{N-3}\ran \] ds
+\dbI_{N-2}  \Big\},\label{cost-(N-3)}
\end{align}
where $X^\dag_\D(\cd)$ is the solution of \rf{state-D-(N-2)} with the initial state $X^\dag_\D(t_{N-2})=X_{N-3}(t_{N-2})$.
Note that
\begin{align}
\dbI_{N-2}&=\dbE\Big\{\int_{t_{N-2}}^{t_{N-1}} \[\lan Q X^\dag_{N-2},X_{N-2}^\dag\ran+\lan\bar Q\dbE_{t_{N-3}}[X_{N-2}^\dag],\dbE_{t_{N-3}}[X_{N-2}^\dag]\ran
 \nn\\
&\qq+\lan\ti Q\dbE[X_{N-2}^\dag],\dbE[X_{N-2}^\dag]\ran+\lan Ru^\dag_{N-2},u^\dag_{N-2}\ran \] ds+\int_{t_{N-1}}^{T} \[\lan Q X^\dag_{N-1},X_{N-1}^\dag\ran
 \nn\\
&\qq+\lan\bar Q\dbE_{t_{N-3}}[X_{N-1}^\dag],\dbE_{t_{N-3}}[X_{N-1}^\dag]\ran+\lan\ti Q\dbE[X_{N-1}^\dag],\dbE[X_{N-1}^\dag]\ran+\lan Ru^\dag_{N-1},u^\dag_{N-1}\ran \] ds\nn\\
&\qq+\lan GX^\dag_{N-1}(T),X^\dag_{N-1}(T)\ran+\lan\bar G\dbE_{t_{N-3}}[X^\dag_{N-1}(T)],\dbE_{t_{N-3}}[X^\dag_{N-1}(T)]\ran\nn\\
&\qq+\lan\ti G\dbE[X^\dag_{N-1}(T)],\dbE[X^\dag_{N-1}(T)]\ran  \Big\}\nn\\
&=\dbE\Big\{\int_{t_{N-2}}^{t_{N-1}} \[\lan Q X^\dag_{N-2},X_{N-2}^\dag\ran+\lan\bar Q\dbE_{t_{N-3}}[X_{N-2}^\dag],\dbE_{t_{N-3}}[X_{N-2}^\dag]\ran
 \nn\\
&\qq+\lan\ti Q\dbE[X_{N-2}^\dag],\dbE[X_{N-2}^\dag]\ran+\lan Ru^\dag_{N-2},u^\dag_{N-2}\ran \] ds\nn\\
&\qq+\big\lan \G_{N-2}(t_{N-1})X^\dag_{N-2}(t_{N-1}),\,X^\dag_{N-2}(t_{N-1})\big\ran\nn\\
&\qq+\big\lan [\bar\G_{N-2}(t_{N-1})-\G_{N-2}(t_{N-1})]\dbE_{t_{N-3}}[X^\dag_{N-2}(t_{N-1})],\,\dbE_{t_{N-3}}
 [X^\dag_{N-2}(t_{N-1})]\big\ran\nn\\
&\qq+\big\lan [\ti\G_{N-2}(t_{N-1})-\bar  \G_{N-2}(t_{N-1})]\dbE[X^\dag_{N-2}(t_{N-1})],\,\dbE[X^\dag_{N-2}(t_{N-1})]\big\ran\Big\}.\nn
\end{align}
Then by Lemma \ref{lem:rep} again, we have
\begin{align}
\dbI_{N-2}&=\dbE\Big\{\big\lan \G_{N-3}(t_{N-2})X_{N-3}(t_{N-2}),\,X_{N-3}(t_{N-2})\big\ran\nn\\
&\qq+\big\lan [\bar\G_{N-3}(t_{N-2})-\G_{N-3}(t_{N-2})]\dbE_{t_{N-3}}[X_{N-3}(t_{N-2})],\,\dbE_{t_{N-3}}
 [X_{N-3}(t_{N-2})]\big\ran\nn\\
&\qq+\big\lan [\ti\G_{N-3}(t_{N-2})-\bar\G_{N-3}(t_{N-2})]
\dbE[X_{N-3}(t_{N-2})],\,\dbE[X_{N-3}(t_{N-2})]\big\ran\Big\},\label{repre-(N-2)}
\end{align}
with the Lyapunov equations over $[t_{N-2},t_{N-1}]$:
\begin{equation*}\left\{\begin{aligned}
   &\dot{\check\G}_{N-3}+\check\G_{N-3} (A+B\Psi_{N-2}^\dag)+(A+B\Psi_{N-2}^\dag)^\top \check\G_{N-3}\\
   &\q+ (C+D\Psi_{N-2}^\dag)^\top \check\G_{N-3}  (C+D\Psi_{N-2}^\dag)+Q+\Psi_{N-2}^{^\dag\top} R\Psi_{N-2}^\dag=0, \\
   &\check\G_{N-3}(t_{N-1})=\G_{N-2}(t_{N-1}),
\end{aligned}\right.\end{equation*}
\begin{equation*}\left\{\begin{aligned}
  &\dot{\G}_{N-3}+\G_{N-3} [A+\bar A+B(\Psi_{N-2}^\dag+\bar\Psi_{N-2}^\dag)]\\
   &\q+ [A+\bar A+B(\Psi_{N-2}^\dag+\bar\Psi_{N-2}^\dag)]^\top \G_{N-3}\\
   &\q+ [C+\bar C+D(\Psi_{N-2}^\dag+\bar\Psi_{N-2}^\dag)]^\top\check\G_{N-3} [C+\bar C+D(\Psi_{N-2}^\dag+\bar\Psi_{N-2}^\dag)]\\
   &\q+Q+(\Psi_{N-2}^\dag+\bar\Psi_{N-2}^\dag)^\top R(\Psi_{N-2}^\dag+\bar\Psi_{N-2}^\dag)=0, \\
   &\G_{N-3}(t_{N-1})=\G_{N-2}(t_{N-1}),
\end{aligned}\right.\end{equation*}
\begin{equation*}\left\{\begin{aligned}
   &\dot{\bar\G}_{N-3}+\bar\G_{N-3} [A+\bar A+B(\Psi_{N-2}^\dag+\bar\Psi_{N-2}^\dag)]\\
   &\q+ [A+\bar A+B(\Psi_{N-2}^\dag+\bar\Psi_{N-2}^\dag)]^\top\bar \G_{N-3}\\
   &\q+ [C+\bar C+D(\Psi_{N-2}^\dag+\bar\Psi_{N-2}^\dag)]^\top\check\G_{N-3} [C+\bar C+D(\Psi_{N-2}^\dag+\bar\Psi_{N-2}^\dag)]\\
   &\q+Q+\bar Q+(\Psi_{N-2}^\dag+\bar\Psi_{N-2}^\dag)^\top R(\Psi_{N-2}^\dag+\bar\Psi_{N-2}^\dag)=0, \\
   &\bar\G_{N-3}(t_{N-1})=\bar\G_{N-2}(t_{N-1}),
\end{aligned}\right.\end{equation*}
\begin{equation*}\left\{\begin{aligned}
   &\dot{\ti\G}_{N-3}+\ti\G_{N-3} [A+\bar A+\ti A+B(\Psi_{N-2}^\dag+\bar\Psi_{N-2}^\dag+\ti\Psi_{N-2}^\dag)]\\
   &\q+[A+\bar A+\ti A+B(\Psi_{N-2}^\dag+\bar\Psi_{N-2}^\dag+\ti\Psi_{N-2}^\dag)]^\top\ti \G_{N-3}\\
   &\q+ [C+\bar C+\ti C+D(\Psi_{N-2}^\dag+\bar\Psi_{N-2}^\dag+\ti\Psi_{N-2}^\dag)]^\top\check\G_{N-3}\\
   &\q\times  [C+\bar C+\ti C+D(\Psi_{N-2}^\dag+\bar\Psi_{N-2}^\dag+\ti\Psi_{N-2}^\dag)]+Q+\bar Q+\ti Q\\
   &\q+(\Psi_{N-2}^\dag+\bar\Psi_{N-2}^\dag+\ti\Psi_{N-2}^\dag)^\top R(\Psi_{N-2}^\dag+\bar\Psi_{N-2}^\dag+\ti\Psi_{N-2}^\dag)=0, \\
   &\ti\G_{N-3}(t_{N-1})=\ti\G_{N-2}(t_{N-1}).
\end{aligned}\right.\end{equation*}
By substituting \rf{repre-(N-2)} into \rf{cost-(N-3)}, we get
\begin{align}
&J_{N-3}(t_{N-3},\xi;u(\cd))
=\dbE\Big\{\int_{t_{N-3}}^{t_{N-2}} \[\lan Q X_{N-3},X_{N-3}\ran+\lan\bar Q\dbE_{t_{N-3}}[X_{N-3}],\dbE_{t_{N-3}}[X_{N-3}]\ran  \nn\\
&\q+\lan\ti Q\dbE[X_{N-3}],\dbE[X_{N-3}]\ran+\lan Ru_{N-3},u_{N-3}\ran \] ds+\big\lan \G_{N-3}(t_{N-2})X_{N-3}(t_{N-2}),\,X_{N-3}(t_{N-2})\big\ran\nn\\
&\q+\big\lan [\bar\G_{N-3}(t_{N-2})-\G_{N-3}(t_{N-2})]\dbE_{t_{N-3}}[X_{N-3}(t_{N-2})],\,\dbE_{t_{N-3}}
 [X_{N-3}(t_{N-2})]\big\ran\nn\\
&\q+\big\lan [\ti\G_{N-3}(t_{N-2})-\bar\G_{N-3}(t_{N-2})]
\dbE[X_{N-3}(t_{N-2})],\,\dbE[X_{N-3}(t_{N-2})]\big\ran  \Big\}.\label{cost-(N-3)-1}
\end{align}
With the state equation \rf{state-(N-3)} and the cost functional \rf{cost-(N-3)-1},
we get a classical mean-field SLQ optimal control over $[t_{N-3},t_{N-2}]$.
By Theorem \ref{Prop:pre-solution-clsoed} again,  the unique optimal control
$ u^\dag_{N-3}(\cd)$ on $[t_{N-3},t_{N-2})$ is given by
\begin{align*}
u^\dag_{N-3} (s)&=\Psi^\dag_{N-3}(s)X^\dag_{N-3}(s)+\bar\Psi^\dag_{N-3}(s)\dbE_{t_{N-3}}[X^\dag_{N-3}(s)]\\
&\q+\ti \Psi^\dag_{N-3}(s)\dbE[X^\dag_{N-3}(s)],  \q s\in[t_{N-3},t_{N-2}),
\end{align*}
where
\begin{equation*}\left\{\begin{aligned}
  & dX^\dag_{N-3}(s) =\big\{(A+B\Psi^\dag_{N-3})X^\dag_{N-3}+ (\bar A+B\bar\Psi^\dag_{N-3})\dbE_{t_{N-3}}[X^\dag_{N-3}]\\
   &\qq\qq\q~+ (\ti A+B\ti\Psi^\dag_{N-3})\dbE[X^\dag_{N-3}]\big\}ds\\
   &\qq\qq\q~+\big\{(C+D\Psi^\dag_{N-3})X^\dag_{N-3}+ (\bar C+D\bar\Psi^\dag_{N-3})\dbE_{t_{N-3}}[X^\dag_{N-3}]\\
   &\qq\qq\q~+ (\ti C+D\ti\Psi^\dag_{N-3})\dbE[X^\dag_{N-3}]\big\}dW(s), \q s\in[t_{N-3},t_{N-2}),\\
    &X^\dag_{N-3}(t_{N-3}) = \xi,
\end{aligned}\right.\end{equation*}
and
\begin{align*}
\Psi_{N-3}^\dag&=-\big[R+D^\top P_{N-3}D\big]^{-1}\big[B^\top P_{N-3}+D^\top P_{N-3}C\big],\\
\bar\Psi_{N-3}^\dag&=-\big[R+D^\top P_{N-3}D\big]^{-1}\big[B^\top(\Pi_{N-3}- P_{N-3})+D^\top P_{N-3}\bar C\big],\\
\ti\Psi_{N-3}^\dag&=-\big[R+D^\top P_{N-3}D\big]^{-1}\big[B^\top(\Phi_{N-3}- \Pi_{N-3})+D^\top P_{N-3}\ti C\big],
\end{align*}
with the Riccati equations $[t_{N-3},t_{N-2}]$:
\begin{equation*}\left\{\begin{aligned}
   &\dot{P}_{N-3}+P_{N-3}A+A^\top P_{N-3}+C^\top P_{N-3}C+Q\\
   &\q-(P_{N-3}B+C^\top P_{N-3}D)(R+D^\top P_{N-3}D)^{-1}(B^\top P_{N-3}+D^\top P_{N-3}C)=0, \\
   &P_{N-3}(t_{N-2})=\G_{N-3}(t_{N-2}),
\end{aligned}\right.\end{equation*}
\begin{equation*}\left\{\begin{aligned}
   &\dot{\Pi}_{N-3}+\Pi_{N-3}(A+\bar A)+(A+\bar A)^\top \Pi_{N-3}+(C+\bar C)^\top P_{N-3}(C+\bar C)\\
   &\q+Q+\bar Q-[\Pi_{N-3} B+(C+\bar C)^\top P_{N-3}D](R+D^\top P_{N-3}D)^{-1}\\
   &\q\times[B^\top \Pi_{N-3}+D^\top P_{N-3}(C+\bar C)]=0,\\
   &\Pi_{N-3}(t_{N-2})=\bar\G_{N-3}(t_{N-2}),
\end{aligned}\right.\end{equation*}
\begin{equation*}\left\{\begin{aligned}
  &\dot{\Phi}_{N-3}+\Phi_{N-3}(A+\bar A+\ti A)+(A+\bar A+\ti A)^\top \Phi_{N-3}\\
   &\qq+(C+\bar C+\ti C)^\top P_{N-3}(C+\bar C+\ti C)+Q+\bar Q+\ti Q\\
   &\qq-[\Phi_{N-3} B+(C+\bar C+\ti C)^\top P_{N-3}D](R+D^\top P_{N-3}D)^{-1}\\
   &\qq\times[B^\top \Phi_{N-3}+D^\top P_{N-3}(C+\bar C+\ti C)]=0,\\
   &\Phi_{N-3}(t_{N-2})=\ti\G_{N-3}(t_{N-2}).
\end{aligned}\right.\end{equation*}
For $\f(\cd)=P(\cd),\Pi(\cd),\Phi(\cd),X^\dag(\cd),u^\dag(\cd),\Psi^\dag(\cd),\bar\Psi^\dag(\cd),\ti\Psi^\dag(\cd)$,
we extend $\f_\D(\cd)$ from $[t_{N-2},T]$ to $[t_{N-3},T]$ by the following:
\begin{equation*}
\f_\D(s)=\left\{\begin{aligned}
  &\f_\D(s),\qq s\in(t_{N-2},T];\\
   &\f_{N-3}(s),\q s\in[t_{N-3},t_{N-2}].
\end{aligned}\right.
\end{equation*}
For $\a(\cd)=\check\G(\cd),\G(\cd),\bar\G(\cd),\ti\G(\cd)$,
we define:
\begin{equation*}
\a_\D(s)=\left\{\begin{aligned}
  &\a_{N-2}(s),\q s\in(t_{N-1},T];\\
   &\a_{N-3}(s),\q s\in[t_{N-2},t_{N-1}].
\end{aligned}\right.
\end{equation*}

\ms

\textbf{Problem on $[0,T]$.}
By continuing the above process,  we have the control process:
\begin{align*}
u^\dag_\D (s)&=\Psi^\dag_\D(s)X^\dag_\D(s)+\bar\Psi^\dag_\D(s)\dbE_{\rho_\D(s)}[X^\dag_\D(s)]+\ti \Psi^\dag_\D(s)\dbE[X^\dag_\D(s)],  \q s\in[0,T],
\end{align*}
and the state equation:
\bel{state-D}\left\{\begin{aligned}
  & dX^\dag_\D(s) =\big\{(A+B\Psi^\dag_\D)X^\dag_\D+ (\bar A+B\bar\Psi^\dag_\D)\dbE_{\rho_\D(s)}[X^\dag_\D]+ (\ti A+B\ti\Psi^\dag_\D)\dbE[X^\dag_\D]\big\}ds\\
   &\qq\qq\q~+\big\{(C+D\Psi^\dag_\D)X^\dag_\D+ (\bar C+D\bar\Psi^\dag_\D)\dbE_{\rho_\D(s)}[X^\dag_\D]\\
   &\qq\qq\q~+ (\ti C+D\ti\Psi^\dag_\D)\dbE[X^\dag_\D]\big\}dW(s), \q s\in[0,T],\\
    &X^\dag_\D(0) = \xi,
\end{aligned}\right.\ee
where the feedback strategies are given by
\begin{align}
\Psi^\dag_\D&=-\big[R+D^\top P_\D D\big]^{-1}\big[B^\top P_\D+D^\top P_\D C],\nn\\
\bar\Psi^\dag_\D&=-\big[R+D^\top P_\D D\big]^{-1}\big[B^\top(\Pi_\D- P_\D)+D^\top P_\D\bar C],\nn\\
\ti\Psi^\dag_\D&=-\big[R+D^\top P_\D D\big]^{-1}\big[B^\top(\Phi_\D- \Pi_\D)+D^\top P_\D \ti C],
\label{Psi-D}
\end{align}
the Lyapunov equations over $[t_1,T]$ are given by
\bel{LyE-1-D}\left\{\begin{aligned}
   &\dot{\check\G}_\D+\check\G_\D (A+B\Psi^\dag_\D)+(A+B\Psi^\dag_\D)^\top \check\G_\D\\
   &\q+ (C+D\Psi^\dag_\D)^\top \check\G_\D  (C+D\Psi^\dag_\D)+Q+\Psi_\D^{^\dag\top} R\Psi^\dag_\D=0, \q s\in \bigcup (t_{k-1},t_k),\\
   &\check\G_\D(t_{k})= \G_\D(t_{k}),
\end{aligned}\right.\ee
\bel{LyE-2-D}\left\{\begin{aligned}
  &\dot{\G}_\D+\G_\D [A+\bar A+B(\Psi^\dag_\D+\bar\Psi^\dag_\D)]+ [A+\bar A+B(\Psi^\dag_\D+\bar\Psi^\dag_\D)]^\top \G_\D\\
   &\q+ [C+\bar C+D(\Psi^\dag_\D+\bar\Psi^\dag_\D)]^\top\check\G_\D [C+\bar C+D(\Psi^\dag_\D+\bar\Psi^\dag_\D)]\\
   &\q+Q+(\Psi^\dag_\D+\bar\Psi^\dag_\D)^\top R(\Psi^\dag_\D+\bar\Psi^\dag_\D)=0, \\
   &\G_\D(T)=G,
\end{aligned}\right.\ee
\bel{LyE-3-D}\left\{\begin{aligned}
   &\dot{\bar\G}_\D+\bar\G_\D [A+\bar A+B(\Psi^\dag_\D+\bar\Psi^\dag_\D)]+ [A+\bar A+B(\Psi^\dag_\D+\bar\Psi^\dag_\D)]^\top\bar \G_\D\\
   &\q+ [C+\bar C+D(\Psi^\dag_\D+\bar\Psi^\dag_\D)]^\top\check\G_\D [C+\bar C+D(\Psi^\dag_\D+\bar\Psi^\dag_\D)]\\
   &\q+Q+\bar Q+(\Psi^\dag_\D+\bar\Psi^\dag_\D)^\top R(\Psi^\dag_\D+\bar\Psi^\dag_\D)=0, \\
   &\bar\G_\D(T)=G+\bar G,
\end{aligned}\right.\ee
\bel{LyE-4-D}\left\{\begin{aligned}
   &\dot{\ti\G}_\D+\ti\G_\D [A+\bar A+\ti A+B(\Psi^\dag_\D+\bar\Psi^\dag_\D+\ti\Psi^\dag_\D)]\\
   &\q+[A+\bar A+\ti A+B(\Psi^\dag_\D+\bar\Psi^\dag_\D+\ti\Psi^\dag_\D)]^\top\ti \G_\D\\
   &\q+ [C+\bar C+\ti C+D(\Psi^\dag_\D+\bar\Psi^\dag_\D+\ti\Psi^\dag_\D)]^\top\check\G_\D\\
   &\q\times  [C+\bar C+\ti C+D(\Psi^\dag_\D+\bar\Psi^\dag_\D+\ti\Psi^\dag_\D)]+Q+\bar Q+\ti Q\\
   &\q+(\Psi^\dag_\D+\bar\Psi^\dag_\D+\ti\Psi^\dag_\D)^\top R(\Psi^\dag_\D+\bar\Psi^\dag_\D+\ti\Psi^\dag_\D)=0, \\
   &\ti\G_\D(T)=G+\bar G+\ti G,
\end{aligned}\right.\ee
and the Riccati equations over $[0,T]$ are given by
\bel{RE-D1}\left\{\begin{aligned}
   &\dot{P}_\D+P_\D A+A^\top P_\D+C^\top P_\D C+Q-(P_\D B+C^\top P_\D D)\\
   &\q\times(R+D^\top P_\D D)^{-1}(B^\top P_\D+D^\top P_\D C)=0,\q s\in \bigcup (t_{k-1},t_k), \\
   &P_\D (t_k)=\G_\D(t_k),
\end{aligned}\right.\ee
\bel{RE-D2}\left\{\begin{aligned}
   &\dot{\Pi}_\D+\Pi_\D(A+\bar A)+(A+\bar A)^\top \Pi_\D+(C+\bar C)^\top P_\D(C+\bar C)\\
   &\q+Q+\bar Q-[\Pi_\D B+(C+\bar C)^\top P_\D D](R+D^\top P_\D D)^{-1}\\
   &\q\times[B^\top \Pi_\D+D^\top P_\D(C+\bar C)]=0,\q s\in \bigcup (t_{k-1},t_k),\\
   &\Pi_\D (t_k)=\bar\G_\D (t_k),
\end{aligned}\right.\ee
\bel{RE-D3}\left\{\begin{aligned}
  &\dot{\Phi}_\D+\Phi_\D(A+\bar A+\ti A)+(A+\bar A+\ti A)^\top \Phi_\D\\
   &\qq+(C+\bar C+\ti C)^\top P_\D(C+\bar C+\ti C)+Q+\bar Q+\ti Q\\
   &\qq-[\Phi_\D B+(C+\bar C+\ti C)^\top P_\D D](R+D^\top P_\D D)^{-1}\\
   &\qq\times[B^\top \Phi_\D+D^\top P_\D (C+\bar C+\ti C)]=0,\q s\in \bigcup (t_{k-1},t_k),\\
   &\Phi_\D (t_k)=\ti\G_\D (t_k).
\end{aligned}\right.\ee

\ms

\textbf{Convergence.}
Finally, we study the convergence of the functions constructed in the above.
Noting from \rf{LyE-1-D} and \rf{RE-D1} that
\begin{align*}
   &\dot{P}_\D+P_\D A+A^\top P_\D+C^\top P_\D C+Q-(P_\D B+C^\top P_\D D)\\
   &\times(R+D^\top P_\D D)^{-1}(B^\top P_\D+D^\top P_\D C)\\
   &\q=\dot{P}_\D+P_\D (A+B\Psi^\dag_\D)+(A+B\Psi^\dag_\D)^\top P_\D+ (C+D\Psi^\dag_\D)^\top P_\D  (C+D\Psi^\dag_\D)\\
   &\qq+Q+\Psi_\D^{^\dag\top} R\Psi^\dag_\D,
\end{align*}
and $P_\D (t_k)=\G_\D(t_k)$, we have
\bel{P-D-G-D}
P_\D(s)=\check\G_\D(s),\q s\in[0,T].
\ee
Moreover, observe that
\begin{align*}
   &\dot{\Pi}_\D+\Pi_\D(A+\bar A)+(A+\bar A)^\top \Pi_\D+(C+\bar C)^\top P_\D(C+\bar C)+Q+\bar Q\\
   &-[\Pi_\D B+(C+\bar C)^\top P_\D D](R+D^\top P_\D D)^{-1}[B^\top \Pi_\D+D^\top P_\D(C+\bar C)]\\
    &\q=\dot\Pi_\D+\Pi_\D [A+\bar A+B(\Psi^\dag_\D+\bar\Psi^\dag_\D)]+ [A+\bar A+B(\Psi^\dag_\D+\bar\Psi^\dag_\D)]^\top\Pi_\D\\
   &\qq+ [C+\bar C+D(\Psi^\dag_\D+\bar\Psi^\dag_\D)]^\top P_\D [C+\bar C+D(\Psi^\dag_\D+\bar\Psi^\dag_\D)]\\
   &\qq+Q+\bar Q+(\Psi^\dag_\D+\bar\Psi^\dag_\D)^\top R(\Psi^\dag_\D+\bar\Psi^\dag_\D)\\
   &\q=\dot\Pi_\D+\Pi_\D [A+\bar A+B(\Psi^\dag_\D+\bar\Psi^\dag_\D)]+ [A+\bar A+B(\Psi^\dag_\D+\bar\Psi^\dag_\D)]^\top\Pi_\D\\
   &\qq+ [C+\bar C+D(\Psi^\dag_\D+\bar\Psi^\dag_\D)]^\top \check\G_\D [C+\bar C+D(\Psi^\dag_\D+\bar\Psi^\dag_\D)]\\
   &\qq+Q+\bar Q+(\Psi^\dag_\D+\bar\Psi^\dag_\D)^\top R(\Psi^\dag_\D+\bar\Psi^\dag_\D),
\end{align*}
where the last equality is due to \rf{P-D-G-D}. Then with the boundary condition
$\Pi_\D (t_k)=\bar\G_\D (t_k)$ (see \rf{RE-D2}), by comparing the above with \rf{LyE-3-D}, we have
\bel{Pi-D-G-D}
\Pi_\D(s)=\bar\G_\D(s),\q s\in[0,T].
\ee
By the same arguments as the above, we have
\bel{Phi-D-G-D}
\Phi_\D(s)=\ti\G_\D(s),\q s\in[0,T].
\ee
With \rf{P-D-G-D}--\rf{Phi-D-G-D},  we can rewrite \rf{Psi-D} as follows:
\begin{align}
\Psi^\dag_\D&=-\big[R+D^\top \check\G_\D D\big]^{-1}\big[B^\top\check \G_\D+D^\top\check \G_\D C],\nn\\
\bar\Psi^\dag_\D&=-\big[R+D^\top\check \G_\D D\big]^{-1}\big[B^\top(\bar\G_\D- \check\G_\D)+D^\top \check\G_\D\bar C],\nn\\
\ti\Psi^\dag_\D&=-\big[R+D^\top\check \G_\D D\big]^{-1}\big[B^\top(\ti\G_\D- \bar\G_\D)+D^\top \check\G_\D \ti C].
\label{Psi-D-G-re}
\end{align}

Next, we introduce the following equations:
\bel{RE-Xi}\left\{\begin{aligned}
   &\dot{\Xi}_\D+\Xi_\D A+A^\top \Xi_\D+C^\top \Xi_\D C+Q=0,\q s\in \bigcup (t_{k-1},t_k), \\
   &\Xi_\D (t_k)=\G_\D(t_k),
\end{aligned}\right.\ee
\bel{RE-bar-Xi}\left\{\begin{aligned}
   &\dot{\bar\Xi}_\D+\bar\Xi_\D(A+\bar A)+(A+\bar A)^\top \bar\Xi_\D+(C+\bar C)^\top \Xi_\D(C+\bar C)+Q+\bar Q=0,\\
   &\bar\Xi_\D (T)=G+\bar G.
\end{aligned}\right.\ee
By the comparison theorem of Lyapunov equations (see \cite[Proposition 3.2]{Yong2017}, for example),
on $(t_{N-1},t_N]$, we have
$$
0\les P_\D(s)\les \Xi_\D (s), \q 0\les \Pi_\D(s)\les \bar\Xi_\D (s).
$$
At the point $t_{N-1}$,
$$
P_\D(t_{N-1})=\G_\D(t_{N-1}),\q
0\les \Pi_\D(t_{N-1})=\bar\G_\D (t_{N-1})\les \bar\Xi_\D (t_{N-1}).
$$
Thus,
$$
0\les P_\D(s)\les \Xi_\D (s), \q 0\les \Pi_\D(s)\les \bar\Xi_\D (s),\q s\in(t_{N-2},t_{N-1}].
$$
By induction, we have
\bel{G-barG-est}
0\les P_\D(s)=\check\G_\D(s)\les \Xi_\D (s), \q 0\les \Pi_\D(s)=\bar\G_\D(s)\les \bar\Xi_\D (s),\q s\in[0,T].
\ee
Next, by comparing \rf{LyE-2-D} with \rf{LyE-3-D} (noting that $\bar Q(\cd)\ges 0$ and $\bar G\ges 0$), we have
\bel{checkG-est1}
 0\les \G_\D(s)\les \bar\G_\D(s),\q s\in[0,T],
\ee
and then
\bel{checkG-est}
 0\les \G_\D(s)\les \bar\Xi_\D (s),\q s\in[0,T].
\ee
From \rf{RE-Xi} and \rf{checkG-est}, we have
$$
|\Xi_\D(s)|\les K(\|\D\|+|\G_\D(t_k)|)\les K(\|\D\|+|\bar\Xi_\D(t_k)|),\q  s\in(t_{k-1},t_k].
$$
Thus,
\bel{XiD-est}
|\Xi_\D(s)|\les K\Big(1+\sup_{r\in[s,T]}|\bar\Xi_\D(r)|\Big),\q  s\in[0,T].
\ee
From \rf{RE-bar-Xi}, we have
$$
|\bar\Xi_\D(s)|\les K\Big(1+\int_s^T|\Xi_\D(r)|dr\Big)\les K\Big(1+\int_s^T\sup_{\t\in[r,T]}|\bar\Xi_\D(\t)|dr\Big),\q s\in[0,T].
$$
Then by Gr\"{o}nwall inequality, we get
$$
|\bar\Xi_\D(s)|\les K,\q s\in[0,T].
$$
Combining the above with \rf{G-barG-est}--\rf{checkG-est}, we get
$$
|\G_\D(s)|\les |\bar\G_\D(s)|=| \Pi_\D(s)|\les|\bar\Xi_\D(s)|\les K,\q s\in[0,T].
$$
Moreover, by \rf{G-barG-est} and \rf{XiD-est}, we have
$$
 |P_\D(s)|=|\check\G_\D(s)|\les|\Xi_\D(s)| \les K\Big(1+\sup_{r\in[s,T]}|\bar\Xi_\D(r)|\Big)\les K,\q s\in[0,T].
$$
Next, we introduce the following equation:
$$\left\{\begin{aligned}
   &\dot{\ti\Xi}_\D+\ti\Xi_\D(A+\bar A+\ti A)+(A+\bar A+\ti A)^\top \ti\Xi_\D\\
   &\q+(C+\bar C+\ti C)^\top \Xi_\D(C+\bar C+\ti C)+Q+\bar Q+\ti Q=0,\\
   &\bar\Xi_\D (T)=G+\bar G+\ti G.
\end{aligned}\right.$$
By the same argument as the above, we have
$$
|\Phi_\D(s)|=|\ti\G_\D(s)|\les |\ti\Xi_\D(s)|\les K,\q s\in[0,T].
$$
Thus, $P_\D(\cd)$, $\Pi_\D(\cd)$, $\Phi_\D(\cd)$,
$\G_\D(\cd)$, $\check\G_\D(\cd)$, $\bar\G_\D(\cd)$, and $\ti\G_\D(\cd)$ are all uniformly bounded.
It follows that $\Psi_\D(\cd)$,  $\bar\Psi_\D(\cd)$, and $\ti\Psi_\D(\cd)$ are  uniformly bounded, too.
Then from \rf{LyE-2-D}--\rf{LyE-4-D}, we know that
the derivatives of $\G_\D(\cd)$, $\bar\G_\D(\cd)$, and $\ti\G_\D(\cd)$ are uniformly bounded,
and then $\G_\D(\cd)$, $\bar\G_\D(\cd)$, and $\ti\G_\D(\cd)$ are equicontinuous functions.
Then by  Arzel\`{a}-Ascoli theorem, we know that there exist three functions $(\G(\cd),\bar\G(\cd),\ti \G(\cd))\in C([0,T];\dbS^n_+)^3$
such that
$$
\lim_{\|\D\|\to 0}\sup_{s\in[0,T]}\big(|\G_\D(s)-\G(s)|+ |\bar\G_\D(s)-\bar\G(s)|+|\ti\G_\D(s)-\ti\G(s)|\big)=0.
$$
From \rf{LyE-1-D}, we see that $\check\G_\D(t_{k})=\G_\D(t_{k})$ for $k=1,2,...,N$. Note that
the derivatives of $\G_\D(\cd)$ and $\check\G_\D(\cd)$ are uniformly bounded. Then,
$$
|\G_\D(s)- \check\G_\D(s)|\les K\|\D\|,\q s\in[0,T].
$$
Thus,
$$
\lim_{\|\D\|\to 0}\sup_{s\in[0,T]}\big(|\check\G_\D(s)-\G(s)|\big)=0.
$$
Taking $\|\D\|\to 0$ in \rf{LyE-1-D}--\rf{LyE-4-D} and \rf{Psi-D-G-re}, we get that $(\G(\cd),\bar\G(\cd),\ti\G(\cd))$
satisfies \rf{ER-1}--\rf{ER-3}, with $(\Psi^\dag(\cd),\bar\Psi^\dag(\cd),\ti \Psi^\dag(\cd) )$ given by \rf{ES}.
The uniqueness of solutions to \rf{ER-1}--\rf{ER-3} in  $ C([0,T];\dbS^n_+)^3$ can be obtained by a standard method.

\section{Verification Theorem and Local Optimality}\label{sec:verifi}

For any given $t\in[0,T)$ and $\e>0$, consider the following controlled system
\bel{state-e}\left\{\begin{aligned}
   dX^\e(s) &=\big\{A X^\e + \bar A\dbE_t[X^\e]+\ti A\dbE[X^\e]+B u \big\}ds\\
   &\qq+\big\{C X^\e+ \bar C\dbE_t[X^\e]+\ti C\dbE[X^\e]+Du \big\}dW(s), \q s\in[t,t+\e);\\
    dX^\e(s) &=\big\{[A+\bar A] X^\e +\ti A\dbE[X^\e]+B[\Psi^\dag+\bar\Psi^\dag]X^\e +B\ti\Psi^\dag\dbE[X^\e] \big\}ds+\big\{[C + \bar C]X^\e\\
   &\qq+\ti C\dbE[X^\e]+D[\Psi^\dag+\bar\Psi^\dag]X^\e +D\ti\Psi^\dag\dbE[X^\e] \big\}dW(s),\q s\in[t+\e,T],\\
    X^\e(t) &= X^\dag(t),
\end{aligned}\right.\ee
and the cost functional
\begin{align}\label{cost-e}
&J(t,X^\dag(t);u^\e(\cd))
= \dbE\Big\{\int_t^{t+\e} \[\lan QX^\e,X^\e\ran +\lan\bar Q\dbE_t[X^\e],\dbE_t[X^\e]\ran+\lan\ti Q\dbE[X^\e],\dbE[X^\e]\ran\nn\\
&\qq +\lan Ru,u\ran \] ds+\int_{t+\e}^T \[\lan QX^\e,X^\e\ran +\lan\bar Q\dbE_t[X^\e],\dbE_t[X^\e]\ran+\lan\ti Q\dbE[X^\e],\dbE[X^\e]\ran\nn\\
&\qq+\lan R[(\Psi^\dag+\bar\Psi^\dag)X^\e+\ti\Psi^\dag\dbE[X^\e]],[(\Psi^\dag+\bar\Psi^\dag)X^\e+\ti\Psi^\dag\dbE[X^\e]]\ran \] ds \nn\\
&\qq+\lan GX^\e(T),X^\e(T)\ran+\lan\bar G\dbE_t[X^\e(T)],\dbE_t[X^\e(T)]\ran+\lan\ti G\dbE[X^\e(T)],\dbE[X^\e(T)]\ran  \Big\}\nn\\
&\q=: \dbE\Big\{\int_t^{t+\e} \[\lan QX^\e,X^\e\ran +\lan\bar Q\dbE_t[X^\e],\dbE_t[X^\e]\ran+\lan\ti Q\dbE[X^\e],\dbE[X^\e]\ran+\lan Ru,u\ran \] ds+\dbI_\e\Big\}.
\end{align}
By Lemma \ref{lem:rep}, we have
\begin{align}
\dbI_\e&=\dbE\[\big\lan \G(t+\e)X^\e(t+\e),\,X^\e(t+\e)\big\ran+\big\lan [\bar\G(t+\e)-\G(t+\e)]\dbE_{t}[X^\e(t+\e)],\,\dbE_{t}
 [X^\e(t+\e)]\big\ran\nn\\
&\qq+\big\lan [\ti\G(t+\e)-\bar  \G(t+\e)]\dbE[X^\e(t+\e)],\,\dbE[X^\e(t+\e)]\big\ran\].\nn
\end{align}
Then,
\begin{align}\label{cost-e1}
&J(t,X^\dag(t);u^\e(\cd))
= \dbE\[\int_t^{t+\e} \(\lan QX^\e,X^\e\ran +\lan\bar Q\dbE_t[X^\e],\dbE_t[X^\e]\ran\nn\\
&\q+\lan\ti Q\dbE[X^\e],\dbE[X^\e]\ran+\lan Ru,u\ran \) ds +\big\lan \G(t+\e)X^\e(t+\e),\,X^\e(t+\e)\big\ran\nn\\
&\q+\big\lan [\bar\G(t+\e)-\G(t+\e)]\dbE_{t}[X^\e(t+\e)],\,\dbE_{t}[X^\e(t+\e)]\big\ran\nn\\
&\q+\big\lan [\ti\G(t+\e)-\bar  \G(t+\e)]\dbE[X^\e(t+\e)],\,\dbE[X^\e(t+\e)]\big\ran\].
\end{align}
Note that
\begin{align}
&\dbE\big[\big\lan \G(t+\e)X^\e(t+\e),\,X^\e(t+\e)\big\ran\big]=\dbE\[\big\lan \G(t)X^\dag(t),\,X^\dag(t)\big\ran+\int_t^{t+\e}\(\big\lan \dot{\G}X^\e,\,X^\e\big\ran\nn\\
&\q
+2\big\lan \G \big[A X^\e + \bar A\dbE_t[X^\e]+\ti A\dbE[X^\e]+B u\big],\,X^\e\big\ran+\big\lan \G \big[C X^\e + \bar C\dbE_t[X^\e]\nn\\
&\q+\ti C\dbE[X^\e]+D u\big],\, \big[C X^\e + \bar C\dbE_t[X^\e]+\ti C\dbE[X^\e]+D u\big]\big\ran\)ds\];\nn\\
&\dbE\big[\big\lan [\bar\G(t+\e)-\G(t+\e)]\dbE_{t}[X^\e(t+\e)],\,\dbE_{t}[X^\e(t+\e)]\big\ran\big]\nn\\
&\q=\dbE\[\big\lan [\bar\G(t)-\G(t)]X^\dag(t),\,X^\dag(t)\big\ran+\int_t^{t+\e}\(\big\lan [\dot{\bar\G}-\dot{\G}]\dbE_{t}[X^\e],\,\dbE_{t}[X^\e]\big\ran\nn\\
&\qq
+2\big\lan[\bar\G-\G] \big[A \dbE_t[X^\e] + \bar A\dbE_t[X^\e]+\ti A\dbE[X^\e]+B \dbE_t[u]\big],\,\dbE_t[X^\e]\big\ran\)ds\];\nn\\
&\dbE\big[\big\lan [\ti\G(t+\e)-\bar\G(t+\e)]\dbE[X^\e(t+\e)],\,\dbE[X^\e(t+\e)]\big\ran\big]\nn\\
&\q=\dbE\[\big\lan [\ti\G(t)-\bar\G(t)]\dbE[X^\dag(t)],\,\dbE[X^\dag(t)]\big\ran+\int_t^{t+\e}\(\big\lan [\dot{\ti\G}-\dot{\bar\G}]\dbE[X^\e],\,\dbE[X^\e]\big\ran\nn\\
&\qq
+2\big\lan[\ti\G-\bar\G] \big[A \dbE[X^\e] + \bar A\dbE[X^\e]+\ti A\dbE[X^\e]+B \dbE[u]\big],\,\dbE[X^\e]\big\ran\)ds\].\nn
\end{align}
Thus,
\begin{align}
&J(t,X^\dag(t);u^\e(\cd))
= \dbE\[\big\lan \bar\G(t)X^\dag(t),\,X^\dag(t)\big\ran+\big\lan [\ti\G(t)-\bar\G(t)]\dbE[X^\dag(t)],\,\dbE[X^\dag(t)]\big\ran\nn\\
&\q+\int_t^{t+\e} \(\lan QX^\e,X^\e\ran +\lan\bar Q\dbE_t[X^\e],\dbE_t[X^\e]\ran+\lan\ti Q\dbE[X^\e],\dbE[X^\e]\ran+\lan Ru,u\ran \) ds\nn\\
&\q+\int_t^{t+\e}\(\big\lan \dot{\G}X^\e,\,X^\e\big\ran+2\big\lan \G \big[A X^\e + \bar A\dbE_t[X^\e]+\ti A\dbE[X^\e]+B u\big],\,X^\e\big\ran\nn\\
&\q
+\big\lan \G \big[C X^\e + \bar C\dbE_t[X^\e]+\ti C\dbE[X^\e]+D u\big],\, \big[C X^\e + \bar C\dbE_t[X^\e]+\ti C\dbE[X^\e]+D u\big]\big\ran\)ds\nn\\
&\q+\int_t^{t+\e}\(\big\lan [\dot{\bar\G}-\dot{\G}]\dbE_{t}[X^\e],\,\dbE_{t}[X^\e]\big\ran+2\big\lan[\bar\G-\G] \big[A \dbE_t[X^\e] + \bar A\dbE_t[X^\e]\nn\\
&\qq
+\ti A\dbE[X^\e]+B \dbE_t[u]\big],\,\dbE_t[X^\e]\big\ran\)ds\nn\\
&\q+\int_t^{t+\e}\(\big\lan [\dot{\ti\G}-\dot{\bar\G}]\dbE[X^\e],\,\dbE[X^\e]\big\ran
+2\big\lan[\ti\G-\bar\G] \big[A \dbE[X^\e] + \bar A\dbE[X^\e]\nn\\
&\qq+\ti A\dbE[X^\e]+B \dbE[u]\big],\,\dbE[X^\e]\big\ran\)ds\].\label{VT-p1}
\end{align}
By the same arguments as the above, we have
\begin{align}
&J(t,X^\dag(t);u^\dag(\cd))
= \dbE\[\big\lan \bar\G(t)X^\dag(t),\,X^\dag(t)\big\ran+\big\lan [\ti\G(t)-\bar\G(t)]\dbE[X^\dag(t)],\,\dbE[X^\dag(t)]\big\ran\nn\\
&\q+\int_t^{t+\e} \(\lan QX^\dag,X^\dag\ran +\lan\bar Q\dbE_t[X^\dag],\dbE_t[X^\dag]\ran+\lan\ti Q\dbE[X^\dag],\dbE[X^\dag]\ran+\lan Ru^\dag,u^\dag\ran \) ds\nn\\
&\q+\int_t^{t+\e}\(\big\lan \dot{\G}X^\dag,\,X^\dag\big\ran+2\big\lan \G \big[A X^\dag + \bar A\dbE_t[X^\dag]+\ti A\dbE[X^\dag]+B u^\dag\big],\,X^\dag\big\ran\nn\\
&\q
+\big\lan \G \big[C X^\dag + \bar C\dbE_t[X^\dag]+\ti C\dbE[X^\dag]+D u^\dag\big],\,
\big[C X^\dag + \bar C\dbE_t[X^\dag]+\ti C\dbE[X^\dag]+D u^\dag\big]\big\ran\)ds\nn\\
&\q+\int_t^{t+\e}\(\big\lan [\dot{\bar\G}-\dot{\G}]\dbE_{t}[X^\dag],\,\dbE_{t}[X^\dag]\big\ran
+2\big\lan[\bar\G-\G] \big[A \dbE_t[X^\dag] + \bar A\dbE_t[X^\dag]\nn\\
&\qq
+\ti A\dbE[X^\dag]+B \dbE_t[u^\dag]\big],\,\dbE_t[X^\dag]\big\ran\)ds\]\nn\\
&\q+\int_t^{t+\e}\(\big\lan [\dot{\ti\G}-\dot{\bar\G}]\dbE[X^\dag],\,\dbE[X^\dag]\big\ran
+2\big\lan[\ti\G-\bar\G] \big[A \dbE[X^\dag] + \bar A\dbE[X^\dag]\nn\\
&\qq+\ti A\dbE[X^\dag]+B \dbE[u^\dag]\big],\,\dbE[X^\dag]\big\ran\)ds\].\label{VT-p2}
\end{align}
Note that
$$
\lim_{\e\to 0}\dbE\[\sup_{s\in[t,t+\e]}|X^\e(s)-X^\dag(t)|^2\]=0.
$$
Then from \rf{VT-p1} and \rf{VT-p2}, we have
\begin{align}
&\liminf_{\e\to 0}{J(t,X^\dag(t);u^\e(\cd))-J(t,X^\dag(t);u^\dag(\cd))\over\e}\nn\\
&\q= \dbE\[ \lan Ru,u\ran+2\lan  u,B^\top\bar\G X^\dag\ran+2\lan \G [C X^\dag + \bar CX^\dag+\ti C\dbE[X^\dag]],Du\ran\nn\\
&\qq\q+ \lan \G Du,Du\ran+2\big\lan[\ti\G-\bar\G]B \dbE[u],\dbE[X^\dag]\big\ran -\lan Ru^\dag,u^\dag\ran
-2\lan  u^\dag,B^\top\bar\G X^\dag\ran\nn\\
&\qq\q-2\lan \G [C X^\dag + \bar CX^\dag+\ti C\dbE[X^\dag]],Du^\dag\ran -\lan \G Du^\dag,Du^\dag\ran-2\big\lan[\ti\G-\bar\G]B \dbE[u^\dag],\dbE[X^\dag]\big\ran\]\nn\\
&\q= \dbE\[ \big\lan [R+D^\top \G D]u,\,u\big\ran+2\big\lan  u,\,[B^\top\bar\G+D^\top\G C+D^\top\G\bar C]X^\dag \nn\\
&\qq\q+[D^\top\G\ti C+B^\top(\ti\G-\bar\G)]\dbE[X^\dag]\big\ran-\big\lan [R+D^\top \G D]u^\dag,\,u^\dag\big\ran\nn\\
&\qq\q-2\big\lan  u^\dag,\,[B^\top\bar\G+D^\top\G C+D^\top\G\bar C]X^\dag+[D^\top\G\ti C+B^\top(\ti\G-\bar\G)]\dbE[X^\dag]\big\ran \].\nn
\end{align}
Noting that
$$
u^\dag=-[R+D^\top \G D]^{-1}\big\{[B^\top\bar\G+D^\top\G C+D^\top\G\bar C]X^\dag+[D^\top\G\ti C+B^\top(\ti\G-\bar\G)]\dbE[X^\dag]\big\},
$$
which is the minimizer of the mapping
$$
u\mapsto \big\{\big\lan [R+D^\top \G D]u,\,u\big\ran+2\big\lan  u,\,[B^\top\bar\G+D^\top\G C+D^\top\G\bar C]X^\dag +[D^\top\G\ti C+B^\top(\ti\G-\bar\G)]\dbE[X^\dag]\big\ran\big\},
$$
we have
$$
\liminf_{\e\to 0}{J(t,X^\dag(t);u^\e(\cd))-J(t,X^\dag(t);u^\dag(\cd))\over\e}\ges 0.
$$
This completes the proof.

\section{Conclusions}

\ms

In this paper, we have studied linear quadratic optimal control problems. The state equation and  cost functional contain the state and control, together with the expectations and conditional expectations of the state. The problem turns out to be time-inconsistent, meaning that the optimal control selected for the given initial pair will not stay optimal thereafter. The main purpose is to indicate that in our framework, the time-inconsistency is due to the fact that the conditional expectation of the state appears in the state equation and cost functional.

\ms

We have not got into the most generality. There are many possible extensions. We mention some of them:

\ms

$\bullet$ The weight matrix valued maps can depend on two time variables. This is the case if we have general discounting in the cost functional (see \cite{Yong2017}).

\ms

$\bullet$ The state equation can contain the expectations and conditional expectations  of control processes (see \cite{Yong2013,Yong2017}).

\ms

$\bullet$ The state equation can contain nonhomogeneous terms and the cost functional contains linear terms. In this case, proper BSDEs are expected to be involved (See \cite{Sun-Yong2020} for the time-consistent problem).

\ms

$\bullet$ Without assuming \ref{ass:A2}, one has the so-called indefinite LQ problems. We could assume the map $u(\cd)\mapsto J(t,x;u(\cd))$ to be uniformly convex (See also \cite{Sun-Yong2020}).

\ms

$\bullet$ All the involved functions are random.

\section*{Appendix}

\subsection*{A1: Proof of Lemma \ref{lem:rep}}

The well-posedness of \rf{lem:rep1}--\rf{lem:rep4} is standard and we are going to prove the representation \rf{lem:rep5} only.
By applying  It\^{o} formula, we have
\begin{align}
&d\dbE\[\big\lan \check\G(X-\dbE_t[X]),(X-\dbE_t[X])\big\ran
+\big\lan \G(\dbE_t[X]-\dbE_\t[X]),(\dbE_t[X]-\dbE_\t[X])\big\ran\nn\\
&\qq+\big\lan \bar\G(\dbE_\t[X]-\dbE[X]),(\dbE_\t[X]-\dbE[X])\big\ran
+\big\lan \ti\G\dbE[X],\dbE[X]\big\ran \]\nn\\
&\q=\dbE\[\big\lan \dot{\check\G}(X-\dbE_t[X]),(X-\dbE_t[X])\big\ran
+\big\lan \check\G\cA(X-\dbE_t[X]),(X-\dbE_t[X])\big\ran\nn\\
&\qq+\big\lan \check\G(X-\dbE_t[X]),\cA(X-\dbE_t[X])\big\ran
+\big\lan \check\G(\cC X+ \bar \cC\dbE_t[X]+\ti\cC\dbE[X]),\nn\\
&\qq\q(\cC X+ \bar \cC\dbE_t[X]+\ti\cC\dbE[X])\big\ran+\big\lan \dot{\G}(\dbE_t[X]-\dbE_\t[X]),(\dbE_t[X]-\dbE_\t[X])\big\ran\nn\\
&\qq+\big\lan \G(\cA+\bar\cA)(\dbE_t[X]-\dbE_\t[X]),(\dbE_t[X]-\dbE_\t[X])\big\ran\nn\\
&\qq+\big\lan \G(\dbE_t[X]-\dbE_\t[X]),(\cA+\bar\cA)(\dbE_t[X]-\dbE_\t[X])\big\ran\nn\\
&\qq+\big\lan\dot{\bar\G}(\dbE_\t[X]-\dbE[X]),(\dbE_\t[X]-\dbE[X])\big\ran\nn\\
&\qq+\big\lan\bar\G(\cA+\bar\cA)(\dbE_\t[X]-\dbE[X]),(\dbE_\t[X]-\dbE[X])\big\ran\nn\\
&\qq+\big\lan\bar\G(\dbE_\t[X]-\dbE[X]),(\cA+\bar\cA)(\dbE_\t[X]-\dbE[X])\big\ran\nn\\
&\qq+\big\lan\dot{\ti\G}\dbE[X],\dbE[X]\big\ran
+\big\lan\ti\G(\cA+\bar\cA+\ti\cA)\dbE[X],\dbE[X]\big\ran\nn\\
&\qq+\big\lan\ti\G\dbE[X],(\cA+\bar\cA+\ti\cA)\dbE[X]\big\ran \]ds.\nn
\end{align}
Substituting \rf{lem:rep1}--\rf{lem:rep4} into the above yields that
\begin{align}
&d\dbE\[\big\lan \check\G(X-\dbE_t[X]),(X-\dbE_t[X])\big\ran
+\big\lan \G(\dbE_t[X]-\dbE_\t[X]),(\dbE_t[X]-\dbE_\t[X])\big\ran\nn\\
&\qq+\big\lan \bar\G(\dbE_\t[X]-\dbE[X]),(\dbE_\t[X]-\dbE[X])\big\ran
+\big\lan \ti\G\dbE[X],\dbE[X]\big\ran \]\nn\\
&\q=-\dbE\[\big\lan ( \cC^\top \check\G\cC+\cQ_1)(X-\dbE_t[X]),(X-\dbE_t[X])\big\ran\nn\\
&\qq-\big\lan \check\G(\cC X+ \bar \cC\dbE_t[X]+\ti\cC\dbE[X]),(\cC X+ \bar \cC\dbE_t[X]+\ti\cC\dbE[X])\big\ran\nn\\
&\qq+\big\lan [(\cC+\bar\cC)^\top \check\G(\cC+\bar\cC)+\cQ_1+\cQ_2](\dbE_t[X]-\dbE_\t[X]),(\dbE_t[X]-\dbE_\t[X])\big\ran\nn\\
&\qq+\big\lan[(\cC+\bar\cC)^\top \check\G(\cC+\bar\cC)+\cQ_1+\cQ_2+\cQ_3](\dbE_\t[X]-\dbE[X]),(\dbE_\t[X]-\dbE[X])\big\ran\nn\\
&\qq+\big\lan[(\cC+\bar\cC+\ti\cC)^\top \check\G(\cC+\bar\cC+\ti\cC)+\cQ_1+\cQ_2+\cQ_3+\cQ_4]\dbE[X],\dbE[X]\big\ran \]ds\nn\\
&\q=-\dbE\[\lan \cQ_1 X,X\ran+\lan \cQ_2\dbE_t[X],\dbE_t[X]\ran+\lan \cQ_3\dbE_\t[X],\dbE_\t[X]\ran+\lan \cQ_4\dbE[X],\dbE [X]\ran\]ds.\nn
\end{align}
Moreover, by the terminal conditions of \rf{lem:rep1}--\rf{lem:rep4}, we have
\begin{align}
&\dbE\[\lan \cG_1 X(T),X(T)\ran+\lan \cG_3\dbE_\t[X(T)],\dbE_\t[X(T)]\ran+\lan \cG_4\dbE[X(T)],\dbE[X(T)]\ran\]\nn\\
&\q=\dbE\[\lan \check\G(T)(X(T)-\dbE_t[X(T)]),(X(T)-\dbE_t[X(T)])\ran\nn\\
&\qq+\lan \G(T)(\dbE_t[X(T)]-\dbE_\t[X(T)]),(\dbE_t[X(T)]-\dbE_\t[X(T)])\ran\nn\\
&\qq+\lan \bar\G(T)(\dbE_\t[X(T)]-\dbE[X(T)]),(\dbE_\t[X(T)]-\dbE[X(T)])\ran\nn\\
&\qq+\lan \ti\G(T)\dbE[X(T)],\dbE[X(T)]\ran \].\nn
\end{align}
Thus,
\begin{align}
&\dbE\[\lan \cG_1 X(T),X(T)\ran+\lan \cG_3\dbE_\t[X(T)],\dbE_\t[X(T)]\ran+\lan \cG_4\dbE[X(T)],\dbE[X(T)]\ran\]\nn\\
&\q=\dbE\[\lan \G(t)(X(t)-\dbE_\t[X(t)]),(X(t)-\dbE_\t[X(t)])\ran\nn\\
&\qq\q+\lan \bar\G(t)(\dbE_\t[X(t)]-\dbE[X(t)]),(\dbE_\t[X(t)]-\dbE[X(t)])\ran+\lan \ti\G(t)\dbE[X(t)],\dbE[X(t)]\ran \nn\\
&\qq\q-\int_t^T\(\lan \cQ_1 X,X\ran+\lan \cQ_2\dbE_t[X],\dbE_t[X]\ran+\lan \cQ_3\dbE_\t[X],\dbE_\t[X]\ran+\lan \cQ_4\dbE[X],\dbE [X]\ran\)ds\],\nn
\end{align}
which completes the proof.

\subsection*{A2: Proof of Proposition \ref{Prop:pre-solution1}}
For any $\lambda\in\dbR$ and $u(\cd)\in\sU[t,T]$, let $X^\lambda(\cd)$ be the state process
corresponding to the control $u^*(\cd)+\lambda u(\cd)$, that is
\bel{state-lambda}\left\{\begin{aligned}
   dX^\lambda(s) &=\big\{AX^\lambda+ \bar A\dbE_t[X^\lambda]+\ti A\dbE[X^\lambda]+B[u^*+\lambda u] \big\}ds\\
   &\qq+\big\{CX^\lambda+ \bar C\dbE_t[X^\lambda]+\ti C\dbE[X^\lambda]+D[u^*+\lambda u] \big\}dW(s), \q s\in[t,T],\\
    X^\lambda(t) &= \xi.
\end{aligned}\right.\ee
Let $X(\cd)$ be the solution of the following equation:
\bel{state-X}\left\{\begin{aligned}
   dX(s) &=\big\{AX+ \bar A\dbE_t[X]+\ti A\dbE[X]+B u \big\}ds\\
   &\qq+\big\{CX+ \bar C\dbE_t[X]+\ti C\dbE[X]+D u \big\}dW(s),\q s\in[t,T], \\
    X(t) &= 0.
\end{aligned}\right.\ee
Then $X^\l(\cd)=X^*(\cd)+\lambda X(\cd)$. It follows that
\begin{align}
0&\les J(t,\xi;u^*(\cd)+\lambda u(\cd))-J(t,\xi;u^*(\cd))\nn\\
&= \lambda^2\dbE\[\int_t^T \(\lan QX,X\ran +\lan\bar Q\dbE_t[X],\dbE_t[X]\ran+\lan\ti Q\dbE[X],\dbE[X]\ran+\lan Ru,u\ran \) ds\nn\\
&\q +\lan GX(T),X(T)\ran +\lan\bar G\dbE_t[X(T)],\dbE_t[X(T)]\ran
+\lan\ti G\dbE[X(T)],\dbE[X(T)]\ran  \]\nn\\
&\q+2\lambda\dbE\[\int_t^T \(\lan QX^*,X\ran +\lan\bar Q\dbE_t[X^*],\dbE_t[X]\ran+\lan\ti Q\dbE[X^*],\dbE[X]\ran+\lan Ru^*,u\ran \) ds\nn\\
&\q +\lan GX^*(T),X(T)\ran +\lan\bar G\dbE_t[X^*(T)],\dbE_t[X(T)]\ran
+\lan\ti G\dbE[X^*(T)],\dbE[X(T)]\ran  \].\label{prop:proof1}
\end{align}
Note that under \ref{ass:A2}, we have
\begin{align*}
&\dbE\[\int_t^T \(\lan QX,X\ran +\lan\bar Q\dbE_t[X],\dbE_t[X]\ran+\lan\ti Q\dbE[X],\dbE[X]\ran+\lan Ru,u\ran \) ds\nn\\
&\q +\lan GX(T),X(T)\ran +\lan\bar G\dbE_t[X(T)],\dbE_t[X(T)]\ran
+\lan\ti G\dbE[X(T)],\dbE[X(T)]\ran  \]\ges 0.
\end{align*}
Thus, \rf{prop:proof1} holds if and only if
\begin{align}
0&=\dbE\[\int_t^T \(\lan QX^*,X\ran +\lan\bar Q\dbE_t[X^*],\dbE_t[X]\ran+\lan\ti Q\dbE[X^*],\dbE[X]\ran+\lan Ru^*,u\ran \) ds\nn\\
&\q +\lan GX^*(T),X(T)\ran +\lan\bar G\dbE_t[X^*(T)],\dbE_t[X(T)]\ran
+\lan\ti G\dbE[X^*(T)],\dbE[X(T)]\ran  \]\nn\\
&=\dbE\[\int_t^T \(\lan (QX^*+\bar Q\dbE_t[X^*]+\ti Q\dbE[X^*]),X\ran+\lan Ru^*,u\ran \) ds\nn\\
&\q +\lan( GX^*(T)+\bar G\dbE_t[X^*(T)]+\ti G\dbE[X^*(T)]),X(T)\ran  \].\label{prop:proof2}
\end{align}
By applying It\^{o} formula to the mapping $s\mapsto\lan X(s),Y(s)\ran$, we have
\begin{align}
&\dbE\big[\lan Y(T),X(T)\ran\big]=\dbE\[\int_t^T\(\big\lan AX+ \bar A\dbE_t[X]+\ti A\dbE[X]+B u,\,Y\big\ran-\big\lan A^\top Y+ \bar A^\top \dbE_t[Y]\nn\\
&\q+\ti A^\top\dbE[Y]+C^\top Z+ \bar C^\top \dbE_t[Z]+\ti C^\top\dbE[Z] +QX^*+ \bar Q \dbE_t[X^*]+\ti Q\dbE[X^*],\, X\big\ran\nn\\
&\q +\big\lan CX+ \bar C\dbE_t[X]+\ti C\dbE[X]+D u,\, Z\big\ran\)ds\]\nn\\
&=\dbE\[\int_t^T\(\big\lan B u,\,Y\big\ran-\big\lan QX^*+ \bar Q \dbE_t[X^*]+\ti Q\dbE[X^*],\, X\big\ran +\big\lan D u,\, Z\big\ran\)ds\]\nn.
\end{align}
Thus,
\begin{align*}
&\dbE\[\int_t^T \(\lan (QX^*+\bar Q\dbE_t[X^*]+\ti Q\dbE[X^*]),X\ran+\lan Ru^*,u\ran \) ds\nn\\
& +\lan( GX^*(T)+\bar G\dbE_t[X^*(T)]+\ti G\dbE[X^*(T)]),X(T)\ran  \]\nn\\
&\q=\dbE\[\int_t^T \big\lan B^\top Y+D^\top Z+R u^*,\,u\big\ran ds\],
\end{align*}
which, together with \rf{prop:proof2}, implies that \rf{Prop:pre-solution1-main1} holds.

\subsection*{A3: Proof of Theorem \ref{Prop:pre-solution-clsoed}}

The existence and uniqueness of solutions  of \rf{RE-pre-1}--\rf{RE-pre-3}
can be obtained directly from \cite{Yong2013}.
We assume that
\begin{align}\label{Prop:pre-solution-clsoed-p1}
Y&=P(X^*-\dbE_t[X^*])+\Pi(\dbE_t[X^*]-\dbE[X^*])+\Phi\dbE[X^*],
\end{align}
for some deterministic and differentiable functions $P(\cd)$, $\Pi(\cd)$, and $\Phi(\cd)$.
Then,
\begin{align*}
Z&=P\big\{CX^*+ \bar C\dbE_t[X^*]+\ti C\dbE[X^*]+Du^* \big\}.
\end{align*}
Next, from \rf{Prop:pre-solution1-main1}, we have
\begin{align*}
&Ru^*+B^\top Y+D^\top P\big\{CX^*+ \bar C\dbE_t[X^*]+\ti C\dbE[X^*]+Du^* \big\}=0,
\end{align*}
which implies that
\begin{align}
u^*&=-(R+D^\top PD)^{-1}\big\{B^\top Y+D^\top PCX^*+ D^\top P\bar C\dbE_t[X^*]+D^\top P\ti C\dbE[X^*] \big\}\nn\\
&=-(R+D^\top PD)^{-1}\big\{(B^\top P+D^\top PC)(X^*-\dbE_t[X^*])+B^\top\Pi(\dbE_t[X^*]-\dbE[X^*])\nn\\
&\qq+B^\top\Phi\dbE[X^*]+D^\top P(C+\bar C)\dbE_t[X^*]+D^\top P\ti C\dbE[X^*]\big\}.\label{Prop:pre-solution-clsoed-p2}
\end{align}
From \rf{Prop:pre-solution-clsoed-p1}, we have
\begin{align*}
&A^\top Y+ \bar A^\top \dbE_t[Y]+\ti A^\top\dbE[Y]+C^\top Z+\bar C^\top \dbE_t[Z]+\ti C^\top\dbE[Z]+QX^*+\bar Q \dbE_t[X^*]+\ti Q\dbE[X^*]\nn\\
&\q=-\dot{P}(X^*-\dbE_t[X^*])-P(AX^*+Bu^*-A\dbE_t[X^*]-B\dbE_t[u^*])\nn\\
&\qq-\dot{\Pi}(\dbE_t[X^*]-\dbE[X^*])-\Pi(A\dbE_t[X^*]+\bar A\dbE_t[X^*]+B\dbE_t[u^*]-A\dbE[X^*]\nn\\
&\qq-\bar A\dbE[X^*]-B\dbE[u^*])-\dot{\Phi}\dbE[X^*]-\Phi(A\dbE[X^*]+\bar A\dbE[X^*]+\ti A\dbE[X^*]+B\dbE[u^*]).
\end{align*}
Then,
\begin{align*}
&A^\top P(X^*-\dbE_t[X^*])+A^\top\Pi(\dbE_t[X^*]-\dbE[X^*])+A^\top\Phi\dbE[X^*]\nn\\
&+ \bar A^\top \dbE_t[Y]+\ti A^\top\dbE[Y]+C^\top PCX^*+ C^\top P\bar C\dbE_t[X^*]+C^\top P\ti C\dbE[X^*]\nn\\
&-C^\top PD(R+D^\top PD)^{-1}\big\{(B^\top P+D^\top PC)(X^*-\dbE_t[X^*])\nn\\
&+B^\top\Pi(\dbE_t[X^*]-\dbE[X^*])+B^\top\Phi\dbE[X^*]+D^\top P(C+\bar C)\dbE_t[X^*]+D^\top P\ti C\dbE[X^*]\big\} \nn\\
&+\bar C^\top \dbE_t[Z]+\ti C^\top\dbE[Z]+QX^*+\bar Q \dbE_t[X^*]+\ti Q\dbE[X^*]\nn\\
&\q=-\dot{P}(X^*-\dbE_t[X^*])-PA(X^*-\dbE_t[X^*])+PB(R+D^\top PD)^{-1}\nn\\
&\qq\times\big\{(B^\top P+D^\top PC)(X^*-\dbE_t[X^*])+B^\top\Pi(\dbE_t[X^*]-\dbE[X^*])\nn\\
&\qq+B^\top\Phi\dbE[X^*]+D^\top P(C+\bar C)\dbE_t[X^*]+D^\top P\ti C\dbE[X^*]\big\}+PB\dbE_t[u^*]\nn\\
&\qq-\dot{\Pi}(\dbE_t[X^*]-\dbE[X^*])-\Pi(A\dbE_t[X^*]+\bar A\dbE_t[X^*]+B\dbE_t[u^*]-A\dbE[X^*]\nn\\
&\qq-\bar A\dbE[X^*]-B\dbE[u^*])-\dot{\Phi}\dbE[X^*]-\Phi(A\dbE[X^*]+\bar A\dbE[X^*]+\ti A\dbE[X^*]+B\dbE[u^*]).
\end{align*}
Let $P(\cd)$ be the unique solution to the Riccati equation \rf{RE-pre-1}, then the above
can be simplified as follows:
\begin{align*}
&A^\top\Pi(\dbE_t[X^*]-\dbE[X^*])+A^\top\Phi\dbE[X^*]+ \bar A^\top \dbE_t[Y]+\ti A^\top\dbE[Y]+ C^\top P(C+\bar C)\dbE_t[X^*]\nn\\
&+C^\top P\ti C\dbE[X^*]-C^\top PD(R+D^\top PD)^{-1}\big\{B^\top\Pi(\dbE_t[X^*]-\dbE[X^*])+B^\top\Phi\dbE[X^*]\nn\\
&+D^\top P(C+\bar C)\dbE_t[X^*]+D^\top P\ti C\dbE[X^*]\big\} +\bar C^\top \dbE_t[Z]+\ti C^\top\dbE[Z]+(Q+\bar Q) \dbE_t[X^*]+\ti Q\dbE[X^*]\nn\\
&\q=PB(R+D^\top PD)^{-1}\big\{B^\top\Pi(\dbE_t[X^*]-\dbE[X^*])+B^\top\Phi\dbE[X^*]+D^\top P(C+\bar C)\dbE_t[X^*]\nn\\
&\qq+D^\top P\ti C\dbE[X^*]\big\}+PB\dbE_t[u^*]-\dot{\Pi}(\dbE_t[X^*]-\dbE[X^*])-\Pi(A\dbE_t[X^*]+\bar A\dbE_t[X^*]\nn\\
&\qq+B\dbE_t[u^*]-A\dbE[X^*]-\bar A\dbE[X^*]-B\dbE[u^*])-\dot{\Phi}\dbE[X^*]\nn\\
&\qq-\Phi(A\dbE[X^*]+\bar A\dbE[X^*]+\ti A\dbE[X^*]+B\dbE[u^*]).
\end{align*}
Note that
\begin{align}
\dbE_t[Y]&=\Pi(\dbE_t[X^*]-\dbE[X^*])+\Phi\dbE[X^*],\nn\\
\dbE_t[Z]&=P(C+\bar C)\dbE_t[X^*]+P\ti C\dbE[X^*]+PD\dbE_t[u^*],\nn \\
\dbE_t[u^*]&=-(R+D^\top PD)^{-1}\big\{B^\top\Pi(\dbE_t[X^*]-\dbE[X^*])+B^\top\Phi\dbE[X^*]\nn\\
&\qq+D^\top P(C+\bar C)\dbE_t[X^*]+D^\top P\ti C\dbE[X^*]\big\}.\nn
\end{align}
Then,
\begin{align*}
&(A+\bar A)^\top\Pi(\dbE_t[X^*]-\dbE[X^*])+(A+\bar A)^\top\Phi\dbE[X^*]+\ti A^\top\dbE[Y]\nn\\
&+ C^\top P(C+\bar C)(\dbE_t[X^*]-\dbE[X^*])+C^\top P(C+\bar C+\ti C)\dbE[X^*]\nn\\
&-C^\top PD(R+D^\top PD)^{-1}\big\{B^\top\Pi(\dbE_t[X^*]-\dbE[X^*])+B^\top\Phi\dbE[X^*]\nn\\
&+D^\top P(C+\bar C)(\dbE_t[X^*]-\dbE[X^*])+D^\top P(C+\bar C+\ti C)\dbE[X^*]\big\} \nn\\
&+\bar C^\top P(C+\bar C)(\dbE_t[X^*]-\dbE[X^*])+\bar C^\top P(C+\bar C+\ti C)\dbE[X^*]\nn\\
&-\bar C^\top PD(R+D^\top PD)^{-1}\big\{B^\top\Pi(\dbE_t[X^*]-\dbE[X^*])+B^\top\Phi\dbE[X^*]\nn\\
&+D^\top P(C+\bar C)(\dbE_t[X^*]-\dbE[X^*])+D^\top P(C+\bar C+\ti C)\dbE[X^*]\big\}+\ti C^\top\dbE[Z]\nn\\
&+(Q+\bar Q) (\dbE_t[X^*]-\dbE[X^*])+(Q+\bar Q+\ti Q)\dbE[X^*]\nn\\
&\q=-\dot{\Pi}(\dbE_t[X^*]-\dbE[X^*])-\Pi(A+\bar A)(\dbE_t[X^*]-\dbE[X^*])\nn\\
&\qq+\Pi B(R+D^\top PD)^{-1}\big\{B^\top\Pi(\dbE_t[X^*]-\dbE[X^*])+B^\top\Phi\dbE[X^*]\nn\\
&\qq+D^\top P(C+\bar C)(\dbE_t[X^*]-\dbE[X^*])+D^\top P(C+\bar C+\ti C)\dbE[X^*]\big\}\nn\\
&\qq+\Pi B\dbE[u^*]-\dot{\Phi}\dbE[X^*]-\Phi(A\dbE[X^*]+\bar A\dbE[X^*]+\ti A\dbE[X^*]+B\dbE[u^*]).
\end{align*}
Let $\Pi(\cd)$ be the unique solution of \rf{RE-pre-2}, then  the above
can be simplified as follows:
\begin{align*}
&(A+\bar A)^\top\Phi\dbE[X^*]+\ti A^\top\dbE[Y]+C^\top P(C+\bar C+\ti C)\dbE[X^*]\nn\\
&-C^\top PD(R+D^\top PD)^{-1}\big\{B^\top\Phi\dbE[X^*]+D^\top P(C+\bar C+\ti C)\dbE[X^*]\big\} \nn\\
&+\bar C^\top P(C+\bar C+\ti C)\dbE[X^*]-\bar C^\top PD(R+D^\top PD)^{-1}\big\{B^\top\Phi\dbE[X^*]\nn\\
&+D^\top P(C+\bar C+\ti C)\dbE[X^*]\big\}+\ti C^\top\dbE[Z]+(Q+\bar Q+\ti Q)\dbE[X^*]\nn\\
&\q=\Pi B(R+D^\top PD)^{-1}\big\{B^\top\Phi\dbE[X^*]+D^\top P(C+\bar C+\ti C)\dbE[X^*]\big\}\nn\\
&\qq+\Pi B\dbE[u^*]-\dot{\Phi}\dbE[X^*]-\Phi(A\dbE[X^*]+\bar A\dbE[X^*]+\ti A\dbE[X^*]+B\dbE[u^*]).
\end{align*}
Note that
\begin{align}
&\dbE[Y]=\Phi\dbE[X^*],\qq
\dbE[Z]=P(C+\bar C+\ti C)\dbE[X^*]+PD\dbE[u^*],\nn \\
&\dbE[u^*]=-(R+D^\top PD)^{-1}\big\{B^\top\Phi+D^\top P(C+\bar C+\ti C)\big\}\dbE[X^*].\nn
\end{align}
Then
\begin{align}
&(A+\bar A)^\top\Phi\dbE[X^*]+\ti A^\top\Phi\dbE[X^*]+C^\top P(C+\bar C+\ti C)\dbE[X^*]\nn\\
&-C^\top PD(R+D^\top PD)^{-1}\big\{B^\top\Phi+D^\top P(C+\bar C+\ti C)\big\}\dbE[X^*]+\bar C^\top P(C+\bar C+\ti C)\dbE[X^*] \nn\\
&-\bar C^\top PD(R+D^\top PD)^{-1}\big\{B^\top\Phi+D^\top P(C+\bar C+\ti C)\big\}\dbE[X^*]+\ti C^\top P(C+\bar C+\ti C)\dbE[X^*]\nn\\
&-\ti C^\top PD(R+D^\top PD)^{-1}\big\{B^\top\Phi+D^\top P(C+\bar C+\ti C)\big\}\dbE[X^*]+(Q+\bar Q+\ti Q)\dbE[X^*]\nn\\
&\q=-\dot{\Phi}\dbE[X^*]-\Phi(A+\bar A+\ti A)\dbE[X^*]+\Phi B(R+D^\top PD)^{-1}\big\{B^\top\Phi+D^\top P(C+\bar C+\ti C)\big\}\dbE[X^*].\nn
\end{align}
Thus, the above holds when $\Phi(\cd)$ is the unique solution of the Riccati equation \rf{RE-pre-3}.
It means that \rf{Prop:pre-solution-clsoed-p1} holds for the unique solution $(P(\cd),\Pi(\cd),\Phi(\cd))$
of \rf{RE-pre-1}--\rf{RE-pre-3}. Then by \rf{Prop:pre-solution-clsoed-p2}, we get
\rf{Prop:pre-solution-clsoed-main1}. This completes the proof.

\end{document}